\documentclass[hidelinks,11pt]{article} 
\usepackage{fullpage}
\usepackage{microtype}
\usepackage{graphicx}
\usepackage{booktabs} 
\usepackage{hyperref}
\usepackage{xcolor}
\usepackage{xcolor}
\usepackage{subcaption}
\usepackage{caption}
\usepackage{wrapfig}
\hypersetup{
    colorlinks=true,
    linkcolor=red,
    citecolor=blue,
    filecolor=magenta,      
    urlcolor=cyan,
   }

\usepackage[utf8]{inputenc} 
\usepackage[T1]{fontenc}    
\usepackage{url}            
\usepackage{amsfonts}       
\usepackage{nicefrac}       
\usepackage{amsmath,amssymb,amsthm}
\usepackage{color}
\usepackage{mathtools}
\usepackage{enumitem}
\usepackage{epstopdf}
\usepackage{tabulary}
\usepackage{natbib}
\usepackage{chngcntr}
\usepackage{bm}
\usepackage[algo2e,linesnumbered,ruled,vlined,resetcount]{algorithm2e}
\RestyleAlgo{ruled}

\SetCommentSty{mycommfont}

\newtheorem{theorem}{Theorem}
\newtheorem{lemma}{Lemma}

\newtheorem{assumption}{Assumption}

\newtheorem{definition}{Definition}
\newtheorem{remark}{Remark}

\newcommand{\norm}[1]{\ensuremath{\left\| #1\right\|}}

\newcommand{\A}{\ensuremath{\mathcal{A}}}

\newcommand{\N}{\ensuremath{\mathcal{N}}}
\newcommand{\R}{\ensuremath{\mathbb{R}}}

\newcommand{\E}{\ensuremath{\mathbb{E}}}
\newcommand{\I}{\ensuremath{\mathcal{I}}}
\newcommand{\X}{\ensuremath{\mathcal{X}}}
\newcommand{\BZ}{\ensuremath{\mathbb{Z}}}

\renewcommand{\O}{\ensuremath{\mathcal{O}}}

\newcommand{\BB}{\ensuremath{\mathbb{B}}}

\newcommand{\CU}{\ensuremath{\mathcal{U}}}
\newcommand{\CO}{\ensuremath{\mathcal{O}}}

\DeclareMathOperator{\Tr}{\mathrm{tr}}

\DeclareMathOperator*{\argmin}{arg\,min}

\title{Online Convex Optimization with Memory and Limited Predictions}
\author{Zhengmiao Wang, Zhi-Wei Liu, Ming Chi, Xiaoling Wang, Housheng Su, Lintao Ye 
\thanks{Z. Wang, ZW. Liu, M. Chi, H. Su and L. Ye are with the School of Artificial Intelligence and Automation, Huazhong University of Science and Technology, Wuhan 430074, China, email: \{wangzhengmiao,zwliu,chiming,shs,yelintao93\}@.hust.edu.cn. X. Wang is with the College of Automation and College of Artificial Intelligence, Nanjing University of Posts and Telecommunications, Nanjing 210023, China, email: xiaolingwang@njupt.edu.cn.} 
}
\begin{document}
\maketitle

\begin{abstract}
This paper addresses an online convex optimization problem where the cost function at each step depends on a history of past decisions (i.e., memory), and the decision maker has access to limited predictions of future cost values within a finite window. The goal is to design an algorithm that minimizes the dynamic regret against the optimal sequence of decisions in hindsight. To this end, we propose a novel predictive algorithm and establish strong theoretical guarantees for its performance. We show that the algorithm's dynamic regret decays exponentially with the length of the prediction window. Our algorithm comprises two general subroutines of independent interest. The first subroutine solves online convex optimization with memory and bandit feedback, achieving a $\sqrt{TV_T}$-dynamic regret, where $V_T$ measures the variation of the optimal decision sequence. The second is a zeroth-order method that attains a linear convergence rate for general convex optimization, matching the best achievable rate of first-order methods. The key to our algorithm is a novel truncated Gaussian smoothing technique when querying the decision points to obtain the predictions. We validate our theoretical results with numerical experiments.
\end{abstract}

\section{Introduction}
Online decision making is a ubiquitous framework that captures a wide range of real-world applications, such as machine learning, control systems, and finance \cite{cesa2006prediction}. In this framework, a decision maker interacts with an {\it unknown} environment by sequentially choosing a decision point $x_t$ at each time step and incurs a cost $f_t(x_t)$ given by some underlying cost function $f_t(\cdot)$, and the goal is to minimize the accumulative cost over a horizon $t=1,\dots,T$. For convex cost functions $f_t(\cdot)$, the problem is known as Online Convex Optimization (OCO), which has been widely studied in the past decades (see, e.g., \cite{hazan2016introduction,shalev2012online} for a comprehensive survey). The OCO framework may be further split into two settings. In the full information setting, the decision maker receives complete knowledge of the cost function $f_t(\cdot)$  after choosing $x_t$ for time step $t$ \cite{zinkevich2003online}; in the bandit information setting, the decision maker only receives $f_t(x_t)$ (and potentially $f_t(\bar{x}_t)$ at some perturbed points $\bar{x}_t$ from $x_t$) after choosing $x_t$ for time step $t$ \cite{flaxman2005online,agarwal2010optimal,saha2011improved}. Since the decision maker needs to choose the decision points in an online manner (i.e., before knowing all the cost functions), regret is a typical performance metric used for algorithms proposed to solve OCO, which compares the accumulative cost incurred by the algorithm to the minimum achievable cost in hindsight \cite{zinkevich2003online}. For the non-stationary environment \cite{zhao2022non}, the minimum achievable cost typically requires time-varying decisions over the horizon $T$, leading to {\it dynamic regret} \cite{yang2016tracking,jadbabaie2015online,zhao2021bandit}.

In many online decision-making problems, the cost incurred at the current time step depends not only on the decision made for the current time step but also on the decisions made in the past. For example, when the underlying environment is characterized by a dynamical system whose state evolves according to the chosen actions, the cost incurred at the current step is governed by the current action and system state, which in turn depends on past actions and system states. Such a framework of OCO (i.e., OCO with memory) has been studied in \cite{arora2012online,anava2015online,zhou2023efficient} and specialized to the online control problem of linear dynamical systems \cite{agarwal2019online,simchowitz2020improper,yan2023online,ye2024online,ye2024learning}. Due to the long-term dependency of the cost on the decision history, i.e., the current decision also influences cost at future steps, there is a line of research that also considers predictions of the cost in a window of future time steps 
\cite{lin2012online,chen2015online,chen2016using,li2020online,li2019control}. We refer to this framework as OCO with memory and predictions. However, the works mentioned above have certain limitations. In \cite{chen2015online,chen2016using}, the cost is assumed to be only a function of decisions in the current step and one step in the past. The classical control algorithm, model predictive control (MPC) \cite{rawlings2012postface,zeilinger2011real,angeli2016theoretical}, considers the receding horizon optimization for
the next $W$ stages. However, similar to \cite{li2019control,li2020online}, the predictions of the full information of the cost functions are assumed to be available in the MPC algorithms.

Motivated by these limitations in the existing works, we study the problem of OCO with memory and limited (i.e., partial) predictions. At each time step, the decision maker is allowed to only query a {\it limited number of points} and obtain potentially {\it noisy predictions} of the cost functions from a finite-length prediction window at these points, and the incurred cost has a {\it long-term dependency} on the past decisions. Focusing on this partial (or bandit) prediction setting is crucial as it reflects many practical applications where the full cost function or its gradient is unavailable. For instance, when interacting with a complex physical system, a proprietary software, or a black-box simulator, it is often impossible or computationally prohibitive to determine the underlying mathematical model \cite{nesterov2017random,balasubramanian2022zeroth}. Instead, the decision-maker can only perform queries by inputting a set or sequence of actions and observing the resulting cost. Such “query-and-observe” interaction precisely corresponds to the bandit feedback model studied in classic OCO problems \cite{flaxman2005online,hazan2014bandit}. We summarize our approaches and contributions.

$\bullet$ {\bf Algorithm Design and Regret Analysis.} In Section~\ref{sec:algorithm design}, we give an algorithm (Algorithm~\ref{alg:two point feedback}) to solve the problem of OCO with memory and limited predictions. Our algorithm extends that in \cite{li2020online} for the full information setting to the partial information setting when making the predictions. We prove in Section~\ref{sec:regret analysis} that the algorithm enjoys a dynamic regret that decays exponentially with the length of the prediction window and matches the dynamic regret of the algorithm in \cite{li2020online}. Our algorithm contains two general subroutines that can be applied to solve wider classes of problems as we describe next. 

$\bullet$ {\bf OCO with Memory and Bandit feedback.} The first subroutine of Algorithm~\ref{alg:two point feedback} can be applied to general OCO with memory and bandit feedback. We show in Theorem~\ref{thm:initial regret bound} our algorithm achieves a $\sqrt{TV_T}$-dynamic regret, where $V_T$ is the variation length of the optimal solution. Our dynamic regret matches that of the algorithm in \cite{zhao2022non} for OCO with memory and full information of cost functions (up to logarithmic factors in $T$). 

$\bullet$ {\bf Faster Convergence Rate of Zeroth-Order Method.} Our core contribution lies in the second subroutine of Algorithm~\ref{alg:two point feedback}, which can be viewed as a zeroth-order (i.e., derivative-free) method for solving general (offline) convex minimization problems. 
For smooth and strongly convex cost functions, our algorithm achieves the same linear convergence rate as the first-order method using gradient descent \cite{nesterov2013introductory,bubeck2015convex}, which improves the convergence rate of the zeroth-order method proposed in \cite{nesterov2017random}. In particular, the convergence rate of our method is independent of the problem horizon length $T$, while applying the method in \cite{nesterov2017random} yields a convergence rate that scales as $\CO(1/T)$. We achieve this speedup by introducing a novel truncated Gaussian smoothing (formally introduced in Definition~\ref{def:truncated gaussian}) when choosing the query points for the function value oracle, i.e., we perturb the decision point $x_t$ at time step $t$ by truncated Gaussian random vectors and query the perturbed points to obtain predictions of the function values. 

{\bf More Related Work.} First, Bandit OCO and zeroth-order methods for convex minimization are naturally related to each other since in both scenarios, a general recursive way to update the action has the form: $x_{t+1}=x_t-\alpha_t g_t$, where $g_t$ is obtained from the function value oracle and can be viewed as an estimate of the true gradient at $x_t$. A major difference is that the performance of algorithms for bandit OCO is characterized by regret while the performance of algorithms for zeroth-order methods is characterized by convergence rate. However, it has been shown that an algorithm for bandit OCO with regret guarantee can be converted to an algorithm for zeroth-order methods with convergence guarantee \cite{shamir2013complexity,shamir2017optimal}. In this work, we propose an algorithm that combines a bandit OCO algorithm and a zeroth-order algorithm to solve the OCO problem with memory and limited predictions. Second, our algorithm design and analysis rely on the use of truncated Gaussian smoothing. The methods of Gaussian smoothing or smoothing over a sphere have been used in bandit OCO or zeroth-order methods \cite{flaxman2005online,nesterov2013introductory,berahas2022theoretical,balasubramanian2022zeroth}. We show that the truncated Gaussian smoothing inherits the desired properties of Gaussian smoothing while also allowing us to prove a faster convergence rate of the zeroth-order method for smooth and strongly convex cost functions. Finally, several extensions of the framework in this paper can be considered, including extending to OCO with unbounded memory (i.e., the current cost depends on all the past decisions) \cite{kumar2024online} and nonsmooth or nonconvex cost functions \cite{lin2022gradient}.


{\bf Notation and terminology.}
The sets of integers and real numbers are denoted as $\mathbb{Z}$ and $\mathbb{R}$, respectively. 
For a set $\mathcal{X}$, $\mathcal{X}^n$ denotes the corresponding Cartesian product $\X\times\cdots\times\X$.  
For a real number $a$, let $\lfloor a \rfloor$ be the largest integer that is smaller than or equal to $a$. 
For a vector $x\in\R^n$, let $\norm{x}=\sqrt{x^{\top}x}$ be its euclidean norm. Let $\norm{P}_F=\sqrt{\Tr(PP^{\top})}$ be the Frobenius norm of $P\in\R^{n\times m}$. 
Let $I_n$ be an $n$-dimensional identity matrix. Given any integer $n\ge1$, $[n]\triangleq\{1,\dots,n\}$. 

\section{Problem Formulation}\label{sec:problem formulation}
We consider an Online Convex Optimization (OCO) framework over $T\in\BZ_{\ge1}$ time steps. At each $t\in[T]\triangleq\{1,\dots,T\}$, a decision maker (i.e., agent) chooses an action $x_t\in\X\subseteq\R^d$ and incurs a cost $f_t(x_{t-h+1},\dots,x_t)$, where $f_t:\R^{dh}\to\R$ is convex and $h\in\BZ_{\ge1}$. In other words, the cost incurred at any time step $t\in\BZ_{\ge1}$ has memory that depends on the current action $x_t$ and the past actions $x_{t-1},x_{t-2},\dots,x_{t-h+1}$. This ``cost with memory'' structure models a broad class of problems where decisions have long-term effects. Common application scenarios include control systems, where a sequence of actions affects future system states due to the dynamics of the system \cite{li2019control}. 
The memory parameter $h$ is an intrinsic characteristic of the problem environment, representing the duration for which past actions remain influential. In practice, its value can be determined by domain knowledge of the specific system or estimated from historical data \cite{zhao2022non,zhou2023efficient}. To proceed, denote 
\begin{equation}\label{eqn:overall obj}
C_T({\bm x}) \triangleq \sum_{t=1}^T f_t(x_{t-h+1:t}),
\end{equation}
where ${\bm x}=(x_1,\dots,x_T)$, $x_{t-h+1:t}=(x_{t-h+1},x_{t-h+2},$ $\dots,x_{t-1},x_{t})$ and assume throughout this paper that $x_t=\bar{x}_0$ for all $t\le 1$ for some initial condition $\bar{x}_0\in\X$. Now, given the complete knowledge of the functions $f_1(\cdot),\dots,f_T(\cdot)$, one can solve
\begin{equation}\label{eqn:offline obj}
{\bm x}^{\star}\in\argmin_{x_1,\dots,x_T\in\X} C_T(\bm{x}),
\end{equation}
where $\bm{x}^{\star}=(x_1^{\star},\dots,x_T^{\star})$. In the online framework described above, since the agent does not have the complete knowledge of $f_1(\cdot),\dots,f_T(\cdot)$ before making the decision $x_t\in\X$ at time step $t\in\BZ_{\ge1}$, solving \eqref{eqn:offline obj} is impossible. The goal is then to design an online algorithm $\A$ for the agent such that the following regret $R(\A)$ is minimized:
\begin{align}
R(\A)\triangleq C_T(\bm{x}^{\A})-C_T(\bm{x}^{\star}),\label{eqn:def of regret}
\end{align}
where $\bm{x}^{\A}=(x_1,\dots,x_T)$ is chosen by the algorithm. 

Even though solving \eqref{eqn:def of regret} requires complete knowledge of $f_1(\cdot),\dots,f_T(\cdot)$ a priori, in practice one may have access to predictions of $f_t(\cdot),\cdots,f_{t+W-1}(\cdot)$ for some $W\in\BZ_{\ge1}$ when choosing $x_t\in\X$ at time step $t\in\BZ_{\ge1}$ \cite{chen2015online,chen2016using,zeilinger2011real,angeli2016theoretical}. However, the available predictions can be inaccurate and scarce (i.e., limited). Unlike existing works as discussed in the introduction, in this paper, we focus on a general framework to capture the inaccuracy and scarcity of the predictions. Specifically, at each time step $t\in\BZ_{\ge1}$, the predictions can be obtained by first running sequences of actions $\{(x_t^j,\dots,x_{t+W-1}^j)\}_{j=1}^p$ on a simulator and then receiving the corresponding costs $\{l_t(x_{t-h+1:t}^j),\dots,l_{t+W-1}(x_{t+W-h:t+W-1}^j)\}_{j=1}^p$, where $l_k:\X^h\to\R$ for all $k\in\{t,\dots,t+W-1\}$. In other words, the simulator is associated with a model $l_t(\cdot)$ that describes the true cost $f_t(\cdot)$ with some model misspecification and can be viewed as a noisy oracle of the function values of $f_t(\cdot)$. Based on the available information, the agent then chooses an action $x_t\in\X$. Note that for each time step $t\in\BZ_{\ge1}$, one may obtain multiple cost prediction trajectories from the simulator by running multiple sequences of actions on the simulator in parallel (or running a single trajectory on multiple simulators). This prediction model, which allows querying future costs for chosen action sequences, is designed to represent practical scenarios where a high-fidelity simulator or a parallel testing environment is available. For example, in robotics, a planner can query a physics simulator to evaluate the cost of different trajectories, or in energy systems, an operator can simulate various generation schedules to predict their impact before commitment. This type of ``query access'' to a value-based oracle corresponds to the widely-studied bandit feedback setting in the classic OCO framework (without memory or predictions) \cite{flaxman2005online,saha2011improved,hazan2014bandit,hazan2016optimal,zhao2021bandit}, since only values of the cost function are available when the agent makes decisions.
Our work applies this well-established and practical model to the specific context of OCO with memory and predictions. To summarize, the agent does the following in each time step $t\ge1$:
\begin{itemize}
\item At the beginning of time step $t\ge1$, run sequences of actions $\{(x_t^j,\dots,x^j_{t+W-1})\}_{j=1}^{p}$ on simulator $l_t(\cdot)$ and receive the predictions $\{(l_t(x^j_{t-h+1:t}),\dots,$ $l_{t+W-1}(x^j_{t+W-h:t+W-1})\}_{j=1}^p$.
\item Compute $x_t$ based on the predictions.
\item Incur the cost $f_t(x_{t-h+1:t})$ and proceed to the next time step $t+1$.\footnote{We assume that the predictions at time step $t$ can be obtained before the algorithm enters the next time step.}
\end{itemize}

\section{Algorithm Design}\label{sec:algorithm design}
To tackle the limited predictions (i.e., bandit information) in our problem setting, we adopt the zeroth-order (i.e., derivative-free) method that has been used for both bandit OCO and (offline) convex optimization problems (e.g., \cite{flaxman2005online,nesterov2017random}). To this end, we extend the algorithm design in \cite{li2019control,lin2012online} for the full information setting to the bandit information setting. The detailed steps are summarized in Algorithm~\ref{alg:two point feedback}. Using similar arguments to those in \cite{li2020online}, one can argue that the overall algorithm can be viewed as using a zeroth-order method for solving \eqref{eqn:offline obj}, which consists of two phases: Initialization (lines~3-7) and zeroth-order method update (lines~8-20). Specifically, let us denote $\bm{x}^j\triangleq(x_1^j,\dots,x_T^j)$ for any $j\in[0,\dots,K-1]$. In the first phase, Algorithm~\ref{alg:two point feedback} uses an OCO subroutine to generate an initial solution sequence $\bm{x}^0=(x_1^0,\dots,x_T^0)$ that is competitive to the optimal solution $\bm{x}^{\star}$ given by \eqref{eqn:offline obj}, i.e., this sequence serves as a high-quality ``warm-start'' for the subsequent steps of Algorithm~\ref{alg:two point feedback}. In the second phase, Algorithm~\ref{alg:two point feedback} employs a zeroth-order method to iteratively compute the update $\bm{x}^{j+1}=\Pi_{(1-\xi^\prime)\X}(\bm{x}^j-\alpha\bm{g}^j)$ for $j=0,\dots,K-1$, where $\Pi_{(1-\xi^\prime)\X}(\cdot)$ denotes the projection onto set $(1-\xi^\prime)\X^T$  and  $\bm{g}^j=(\bm{g}^j_1,\dots,\bm{g}^j_T)\in\R^{dT}$ with $\bm{g}^j_s=\sum_{k=s}^{s+h-1}g_k^j$ for all $s\in[T]$. In particular, the term $\bm{g}^j$ is obtained from $l_1(\cdot),\dots,l_T(\cdot)$ queried at certain points (see lines~6 and 16 of Algorithm~\ref{alg:two point feedback}) which can be viewed as an estimate of the true gradient $\nabla C_T(\bm{x}^j)$. Note that we simply set $f_t(\bm{x})=l_t(\bm{x})=0$ for all $t\le0$, all $t>T$ and all $\bm{x}\in\X^h$ in Algorithm~\ref{alg:two point feedback}, and recall that we have set $x_t^0=\bar{x}_0\in\X$ for all $t\le 1$. We summarize several key features of Algorithm~\ref{alg:two point feedback}.

{\bf Feasibility of Queried and Action Points.} In lines~7 and 17 of Algorithm~\ref{alg:two point feedback}, we project onto smaller sets $(1-\xi)\X$ and $(1-\xi^{\prime})\X$ compared to the original feasible set $\X$, where $\xi,\xi^{\prime}$ are input parameters to the algorithm. Later in Section~\ref{sec:regret analysis}, we will show that via proper choices of input parameters $\delta,\delta^{\prime},\xi,\xi^{\prime}$ to Algorithm~\ref{alg:two point feedback}, projecting onto these smaller sets ensures that both the action points $x_t^K$ (line~21) and the perturbed points $\bar{x}_r^0,\tilde{x}_r^0,\bar{x}^j_{s+h-1},\tilde{x}^j_{s+h-1},\bar{x}_s^{j+1},\tilde{x}_s^{j+1}$ (lines~4, 13, 19) queried for predictions stay within the feasible set $\X$.

\begin{algorithm2e}[h]
\SetNoFillComment
\caption{OCO with Memory and Limited Predictions}
\label{alg:two point feedback}
\KwIn{$\bar{x}_0\in\X,\eta_t,\alpha,\delta,\delta^{\prime}$.}
Pick $u^0\sim\widetilde{\N}(0,I_d)$\\
\For{$t=2-W$ to $T$}{
    \tcc{Initialize $x_{t+W}^0$ using OCO}
    Let $r = t+W-1$.\\   
    Set $\bar{x}_{r}^{0}=x_{r}^{0}+\delta u^{0}$ and $\tilde{x}_{r}^{0}=x_{r}^{0}-\delta u^{0}$.\\
    Receive predictions
    $l_{r}(\bar{x}^0_{r-h+1:r})$ and $l_{r}(\tilde{x}^0_{r-h+1:r})$.\\
    Set $g_r^0=\frac{1}{2\delta}\big[l_r(\bar{x}_{r-h+1:r}^0)-l_s(\tilde{x}_{r-h+1:r}^0)\big]u^0$.\\
    Update $x_{r+1}^0=\mathbf{\Pi}_{(1-\xi)\X}(x_{r}^0-\eta_{r}g_{r}^0)$.\\
    \tcc{Obtain $x_t^K$ using zeroth-order method}
    Let $K=\lfloor\frac{W}{h-1}\rfloor$.\\
    \For{$j=0,\dots,K-1$}{       
        Let $s=t+W-(j+1)(h-1)-(W-(h-1)K)=t+(K-j-1)(h-1)$.\\
        \If{$j=0$}{
             Pick $u_{s}^{j}\overset{i.i.d}{\sim}\widetilde{\N}(0,I_d)$.\\
             Set $\bar{x}_{s+h-1}^{j}=x_{s+h-1}^{j}+\delta^{\prime} u_{s}^{j}$ and $\tilde{x}_{s+h-1}^{j}=x_{s+h-1}^{j}-\delta^{\prime} u_{s}^{j}$.\\
             Receive predictions $l_{s+h-1}(\bar{x}^j_{s:s+h-1})$ and $l_{s+h-1}(\tilde{x}^j_{s:s+h-1})$.\\
        } 
        \For{$k=s,\dots,s+h-1$}{
             Set $g_k^{j}=\frac{1}{2\delta^{\prime}}\big[l_k(\bar{x}_{k-h+1:k}^{j})-l_k(\tilde{x}_{k-h+1:k}^{j})\big]u_s^j$.
        }
        Update $x_s^{j+1}=\mathbf{\Pi}_{(1-\xi^\prime)\X}\big(x_s^{j}-\alpha\sum_{k=s}^{s+h-1}g_k^{j}\big)$.\\
        Pick $u_{s}^{j+1}\overset{i.i.d}\sim\widetilde{\N}(0,I_d)$.\\
        Set $\bar{x}_{s}^{j+1}=x_{s}^{j+1}+\delta^{\prime} u_{s}^{j+1}$ and $\tilde{x}_{s}^{j+1}=x_{s}^{j+1}-\delta^{\prime} u_{s}^{j+1}$.\\
        Receive $l_{s}(\bar{x}^{j+1}_{s-h+1:s})$ and $l_{s}(\tilde{x}^{j+1}_{s-h+1:s})$.\\  
    }
    Play $x_t^{K}$ and incur the cost $f_t(x_{t-h+1:t}^{K})$.
}
\KwOut{$x_1^K,\dots,x_T^K$.}
\end{algorithm2e}

{\bf Time Complexity.} 
The time complexity of Algorithm~\ref{alg:two point feedback} at each time step $t\ge1$ is dominated by the second subroutine, which performs $K = \lfloor\frac{W}{h-1}\rfloor$ update steps (line~17). Since these updates are derivative-free, the complexity is primarily determined by the number of function value queries needed for predictions. We will see in Section~\ref{sec:regret analysis} that these update steps enable the exponential decay in regret with the prediction window size $W$. Finally, the entire algorithm operates without gradient feedback because it is fundamentally designed for the {\it bandit feedback setting}. A core premise of our problem is that the decision-maker cannot access the analytical form of the cost functions or their gradients. Interaction is limited to querying the system with an action and observing the resulting cost value. To adhere to this practical yet challenging information constraint, both subroutines are necessarily designed as derivative-free methods.

{\bf Two-Point Prediction via $l_t(\cdot)$.} Note that in line~5, $\bar{x}^0_{r-h+1:r}=(x^0_{r-h+1},x^0_{r-h+2},\dots,x^0_{r-1},\bar{x}^0_{r})$ and $\tilde{x}^0_{r-h+1:r}=(x^0_{r-h+1},x^0_{r-h+2},\dots,x^0_{r-1},\tilde{x}^0_{r})$, where $\bar{x}^0_{r}$, $\tilde{x}^0_{r}$ are set in line~4. Similar arguments hold for lines~14 and 20. For each call of line~5, line~14 or line~20, Algorithm~\ref{alg:two point feedback} will query two points for $l_t(\cdot)$ to obtain the corresponding predictions. For any $t\in\{2-W,\dots,T\}$, one can check that to obtain $x_t^K$ in line~21 of Algorithm~\ref{alg:two point feedback}, the algorithm needs to use function values of $l_t(\cdot)$ contained in the set $\I_t$ defined as 
\begin{equation}
\I_t\triangleq\big\{l_k(\bar{x}_{k-h+1:k}^j),l_k(\tilde{x}_{k-h+1:k}^j):k\in\{t+(K-j-1)h,\dots,t+W-1\},j\in\{0,\dots,K-1\}\big\}.
\end{equation}
One can also check that all the function values in $\I_t$ are available via the predictions in lines~5, 14 and 20 of Algorithm~\ref{alg:two point feedback}. Moreover, at any $t\in\{2-W,\dots,T\}$, Algorithm~\ref{alg:two point feedback} only needs to maintain in its memory those function values in $\I_t$ (that are needed to obtain $x_t^K$) and remove the other function values from its memory. In other words, the memory of Algorithm~\ref{alg:two point feedback} may be recursively updated through receiving fresh predictions via $l_t(\cdot)$ and removing the function values of $l_t(\cdot)$ that are no longer useful. Finally, we note that for $t=2-W,\dots,T$, the points associated with the predictions in lines~5, 14 and 20 of Algorithm~\ref{alg:two point feedback} correspond to consecutive points on $2(K+1)$ trajectories given by  $\{\bar{x}_1^j,\bar{x}_2^j,\dots\}_{j=0}^{K}$ and $\{\tilde{x}_1^j,\tilde{x}_2^j,\dots\}_{j=0}^{K}$, and thus we do not need to restart the simulator $l_t(\cdot)$ during the course of Algorithm~\ref{alg:two point feedback}.

{\bf Truncated Gaussian Smoothing.} In Algorithm~\ref{alg:two point feedback} (lines~4, 13 and 19), we use a truncated Gaussian random vector denoted as $u_s^j\sim\widetilde{\N}(0,I_d)$ to perturb the point $x_s^j$ to obtain $\bar{x}_s^j$ and $\tilde{x}_s^j$ which are then fed to the simulator $l_t(\cdot)$. Specifically, for any $i\in[d]$, we let the $i$th margin $u_{s,i}^j$ of $u_{s}^j$ to be a truncated Gaussian random variable obtained from the normal distribution such that $u_{s,i}^j$ belongs to the interval $[-(1/2d^2(2h-1))^{1/4},1/(2d^2(2h-1))^{1/4}]$ and the probability density function (pdf) of $u_{s,i}^j$ (i.e., the $i$th margin of $\widetilde{\N}(0,I_d)$) is given by 
\begin{align}\label{eqn:pdf of truncated gaussian}
p(x)=\frac{1}{\sqrt{2\pi}\kappa}e^{-\frac{x^2}{2}}\ \text{and}\ \kappa\ \text{is such that}
 \int_{-1/(2d^2(2h-1))^{\frac{1}{4}}}^{1/(2d^2(2h-1))^{\frac{1}{4}}}p(x)=1.    
\end{align}
It follows that $\lVert u_s^j\rVert\le1/(2(2h-1))^{1/4}\le1$, since $h>1$. Formally, we introduce the following \cite{johnson1995continuous}.
\begin{definition}
\label{def:truncated gaussian}
We use $u\sim\widetilde{\N}(0,I_d)$ to represent an $n$-dimensional truncated Gaussian random vector whose $i$th margin $u_i$ satisfies $u_i\in[-1/(2d^2(2h-1))^{1/4},1/(2d^2(2h-1))^{1/4}]$ and the pdf of $u_i$ is given by \eqref{eqn:pdf of truncated gaussian} for all $i\in[d]$.
\end{definition}
There exist efficient methods for sampling from a truncated Gaussian distribution (see e.g. \cite{robert1995simulation}). As we discussed in the introduction, using the truncated Gaussian allows us to prove the desired regret guarantee of Algorithm~\ref{alg:two point feedback} in the next section.   

\section{Regret Analysis}\label{sec:regret analysis}
We now analyze the regret of Algorithm~\ref{alg:two point feedback}. To proceed, we introduce the following definitions and assumptions, which are standard in the OCO literature \cite{agarwal2010optimal,agarwal2019online,zhao2022non,flaxman2005online}.
\begin{definition}
\label{def:beta-smooth}
A continuous differentiable function $f:\R^d\to\R$ is $\beta$-smooth if 
\begin{equation*}
    \norm{\nabla f(x)-\nabla f(y)}\le\beta\norm{x-y},\ \forall x,y\in\R^d. \\
\end{equation*}
\end{definition}
\begin{definition}
\label{def:mu-strongly convex}
A differentiable function $f:\R^d\to\R$ is $\mu$-strongly convex if 
\begin{equation*}
f(x)-f(y)\le\langle x-y,\nabla f(x)\rangle-\frac{\mu}{2}\norm{x-y}^2,\ \forall x,y\in\R^d.
\end{equation*}
\end{definition}
\begin{definition}
\label{def:lipschitz} A function $f:\R^d\to\R$ is $G$-Lipschitz if 
\begin{equation*}
|f(x)-f(y)|\le G\norm{x-y},\ \forall x,y\in\R^d.
\end{equation*}
\end{definition}

\begin{assumption}
\label{ass:objective functions}
For any $t\in[T]$, $f_t(\cdot)$ is (a) $\beta$-smooth, (b) $\mu$-strongly convex and (c) $G$-Lipschitz.
\end{assumption}

\begin{assumption}
\label{ass:feasible set}
The feasible set $\X\subseteq\R^d$ is convex, contains the origin, and satisfies  $\max_{x,y\in\X}\norm{x-y}\le D$.
\end{assumption}

\begin{remark}\label{remark:feasible small set}
As we discussed in Section~\ref{sec:algorithm design}, Algorithm~\ref{alg:two point feedback} (lines~7 and 17) projects onto smaller sets $(1-\xi)\X$ and $(1-\xi^{\prime})\X$, where $\xi,\xi^{\prime}\in\R_{\ge0}$. By Assumption~\ref{ass:feasible set}, we know that $(1-\xi)\X$ is also convex and contains the origin. Moreover, for any two $y,z\in(1-\xi)\X$, we can write $y=(1-\xi)x_1$ and $z=(1-\xi)x_2$ for some $x_1, x_2 \in \X$, which implies via Assumption~\ref{ass:feasible set} that $\norm{y-z} = \norm{(1-\xi)(x_1-x_2)} = (1-\xi)\norm{x_1-x_2}\le(1-\xi)D$. 
Similar arguments hold for the set $(1-\xi^\prime)\X$. Now, we formally show how to ensure that perturbed used points in lines~4,13,19 of Algorithm~\ref{alg:two point feedback} stay within the set $\X$ as promised in Section~\ref{sec:algorithm design}. To this end, consider a perturbed point, e.g., $\bar{x}=x + \delta u$, where $x$ is kept within the smaller set $(1-\xi)\mathcal{X}$ via projection and $u\sim\tilde{\N}(0,I_d)$ as per Definition~\ref{def:truncated gaussian}. Recall that $\mathcal{X}$ contains the origin by Assumption~\ref{ass:feasible set}, we know that there exists a ball $r\BB=\{x:\norm{x}\le r\}$ with radius $r>0$ such that $r\BB\subseteq \X$, which implies that any point $y \in (1-\xi)\mathcal{X}$ has a distance of at least $\xi r$ to the boundary of $\mathcal{X}$. Recalling that $\norm{u}\le1$ by definition, to guarantee that $\bar{x}=x+\delta u$ remains within $\mathcal{X}$, we can choose $\delta \le \xi r$. 
Finally, it is worth mentioning that in the two phases of Algorithm~\ref{alg:two point feedback}, we use potentially different $\xi,\xi^\prime$ to obtain the sets $(1-\xi)\X,(1-\xi^\prime)\X$, and potentially different $\delta,\delta^{\prime}$ to obtain the perturbed points. Such choices of $\xi,\xi^{\prime},\delta,\delta^{\prime}$ facilitate the subsequent proof of the algorithm's regret upper bounds.
\end{remark}

We will also assume that the model error of $l_t(\cdot)$ is uniformly bounded.
\begin{assumption}
\label{ass:model error}
For any $t\in[T]$, there exists $\varphi_t\in\R_{\ge0}$ such that $|l_t(\bm{x})-f_t(\bm{x})|\le\varphi_t$ for all $\bm{x}\in\X^h$. 
\end{assumption}

\begin{remark}\label{remark:noise model}
Note that the noisy model in Assumption~\ref{ass:model error} is equivalent to assuming that $l_t(\bm{x})=f_t(\bm{x})+\phi_t$ for all $\bm{x}\in\X^h$ with $\phi_t\in\R$ and $\left | \phi _t \right | \le\varphi_t$, where $\phi_t$ can be either deterministic or stochastic. The bounded noise model has also been considered in e.g. \cite{berahas2022theoretical}. The potentially unbounded stochastic noise model has been considered in e.g. \cite{agarwal2011stochastic}, and we leave the extension to such a noise model to future work.
\end{remark}

We define several auxiliary functions that will be useful in our later analysis.
\begin{definition}\label{def:helper functions}
(a) For any $x\in\R^d$ and any $t\in[T]$, define the $h$-unitary function $f_t^c(x)\triangleq f_t(x,\dots,x)$ and its smoothed version $\hat{f}_t^c(x)\triangleq\E_{v}[f_t^c(x+\delta v)]$, where $v\sim\widetilde{\N}(0,I_d)$ and $\delta\in\R_{>0}$.\footnote{Throughout the paper, we use $\E_{x}[\cdot]$ (resp., $\E[\cdot]$) to denote the expectation with respect to a certain random variable (or vector) $x$ (resp., all the randomness in Algorithm~\ref{alg:two point feedback}).}\\
(b) For any $\bm{x}=(x_1,\dots,x_T)\in\R^{dT}$,  $\hat{C}_T(\bm{x})\triangleq\E_{\bm{v}}[C(\bm{x}+\delta^{\prime} \bm{v})]$, where $\bm{v}\sim\widetilde{\N}(0,I_{dT})$ and $\delta^{\prime}\in\R_{>0}$.
\end{definition}

The following result shows that the functions $\hat{f}_t^c(\cdot)$ and $\hat{C}_T(\cdot)$ inherit the properties from $f_t(\cdot)$; the proof is included in Appendix~\ref{app:preliminary proofs}.
\begin{lemma}
\label{lemma:properties of hat f_t^c}
Suppose Assumption~\ref{ass:objective functions} holds for $f_t(\cdot)$ for all $t\in[T]$. Then,  (a) $\hat{f}_t^c(\cdot)$ is $\beta$-smooth and $\mu$-strongly convex; (b) $C_T(\cdot)$ is $\beta h$-smooth, $\mu$-strongly convex and $G\sqrt{Th}$-Lipschitz; and (c) $\hat{C}_T(\cdot)$ is $\beta h$-smooth and $\mu$-strongly convex.
\end{lemma}

More importantly, one can form unbiased gradient estimators of $\hat{f}_t^c(\cdot)$ and $\hat{C}_T(\cdot)$ using only the function values. Formally, we have the following result; the proof can be found in Appendix~\ref{app:preliminary proofs}.

\begin{lemma}
\label{lemma:zeroth-order gradient}
Under Assumption~\ref{ass:objective functions}(a), the following hold:
\noindent(a) For $f_t^c(\cdot)$ and $\hat{f}_t^c(\cdot)$ given by Definition~\ref{def:helper functions}(a), any $x\in\R^d$ and any $\delta \in \R_{>0}$, 
\begin{equation}\label{eqn:two point estimator for hat f_t^c}
\nabla\hat{f}_t^c(x)=\E_v\Big[\frac{f_t^c(x+\delta v)-f_t^c(x-\delta v)}{2\delta}v\Big]
\end{equation}
and 
\begin{equation}\label{eqn:f_t^c and hat f_t^c distance}
|f_t^c(x)-\hat{f}_t^c(x)|\le\frac{\delta^2}{2}\beta d.
\end{equation}
\noindent(b)For $C_T(\cdot)$ and $\hat{C}_T(\cdot)$ given by Definition~\ref{def:helper functions}(b), any $\bm{x}=(x_1,\dots,x_T)\in\R^{dT}$ and any $\delta^\prime \in \R_{>0}$,
\begin{align}
\frac{\partial\hat{C}_T(\bm{x})}{\partial x_s}=&\E_{v_s}\Big[\sum_{k=s}^{s+h-1}\frac{f_k(\bar{x}_{k-h+1:k})-f_k(\tilde{x}_{k-h+1:k})}{2\delta^{\prime}}v_s\Big]
\qquad\forall s\in[T],    \label{eqn:two point estimator for hat C_T}
\end{align}
\begin{equation}\label{eqn:C_T and hat C_T difference}
|C_T(\bm{x})-\hat{C}_T(\bm{x})|\le\frac{\delta^{\prime2}}{2}\beta hTd,
\end{equation}
and 
\begin{equation}\label{eqn:C_T and hat C_T gradient difference}
\norm{\nabla\hat{C}_T(\bm{x})-\nabla C_T(\bm{x})}\le\frac{\delta^{\prime}}{2}\beta h(Td+3)^{3/2},
\end{equation}
where $\bm{v}=(v_1,\dots,v_T)\sim\widetilde{\N}(0,I_{dT})$, $\bar{x}_t=x_t+\delta^{\prime}v_t$ and $\tilde{x}_t=x_t-\delta^{\prime}v_t$.
\end{lemma}

Recall from Eq.~\eqref{eqn:def of regret} that the regret of Algorithm~\ref{alg:two point feedback} is given by $R(\A)=C_T(\bm{x}^{\A})-C_T(\bm{x}^{\star})$, where $\bm{x}^{\A}=(x_1^K,\dots,x_T^K)$ is chosen by the algorithm and $\bm{x}^{\star}=(x_1^{\star},\dots,x_T^{\star})$ is the optimal solution to \eqref{eqn:overall obj}. Since Algorithm~\ref{alg:two point feedback} contains two phases, we will split our analysis accordingly. 

{\bf Initialization Phase Regret.} We first upper bound the contribution of the initialization step in Algorithm~\ref{alg:two point feedback} to $R(\A)$, where we note from the definition of Algorithm~\ref{alg:two point feedback} that $x_1^0,\dots,x_T^0$ are chosen by the initializing phase of the algorithm.

Note that the modification of projecting onto a smaller set $(1-\xi)\mathcal{X}$ ensures the feasibility of query points in Remark~\ref{remark:feasible small set}. To analyze the regret while avoiding dependencies on the variation (also called path length) of the specific optimal solution within the smaller set, we utilize a scaled version of the true optimal sequence as our comparator. Let $\bm{x}^{\star}=(x_1^{\star},\dots,x_T^{\star})$ be the optimal solution to \eqref{eqn:overall obj}. We define a scaled comparator sequence $\hat{\bm{x}}^{\star}=(\hat{x}_1^{\star},\dots,\hat{x}_T^{\star})$ where $\hat{x}_t^{\star} = (1-\xi)x_t^{\star}$. Since $\mathcal{X}$ is convex and contains the origin, $\hat{\bm{x}}^{\star}$ lies within the smaller feasible set $((1-\xi)\mathcal{X})^T$ as discussed in Remark~\ref{remark:feasible small set}. We then decompose the initialization regret as follows:
\begin{align*}
    R(\mathcal{A}_\mathrm{init}) &= C_T(\bm{x}^{0}) - C_T(\bm{x}^{\star})\\
    &= \underbrace{\left( C_T(\bm{x}^{0}) - C_T(\hat{\bm{x}}^{\star}) \right)}_{\text{Regret against scaled comparator}} + \underbrace{\left( C_T(\hat{\bm{x}}^{\star}) - C_T(\bm{x}^{\star}) \right)}_{\text{Approximation Error}}.
\end{align*}
The following theorem bounds the regret in the initialization phase for Algorithm~\ref{alg:two point feedback}, utilizing the fact that the path length of the scaled comparator is bounded by the path length of the true optimum. The detailed proof can be found in Appendix~\ref{app:proofs of main results}.

\begin{theorem}
\label{thm:initial regret bound}
Suppose Assumptions~\ref{ass:objective functions}-\ref{ass:model error} hold. Let $\eta_t=\frac{1}{t\mu}$ in Algorithm~\ref{alg:two point feedback}. Then, the initialization step of Algorithm~\ref{alg:two point feedback} satisfies that
\begin{align}\nonumber
\E[R(\A_\text{init})]
&\le  \Big(\frac{\sqrt{2}(1-\xi)^2D}{\mu}+(1-\xi)Gh^2\Big)V_{T} + T\delta^2\beta d \\\nonumber
& \quad + \Big(\frac{8G^2h^2+hG^2}{2\mu(2h-1)}+\frac{h^3G^2(1-\xi)D}{\delta(2(2h-1))^{1/2}}+\frac{3G^2h^3}{\mu(2(2h-1))^{1/2}}\Big)(1+\ln T) \\\nonumber
& \quad + \frac{2}{\delta^2\mu\sqrt{2h-1}}\sum_{t=1}^T\varphi_t^2  + \Big(\frac{\sqrt{2}h^2G(1-\xi)D}{2\delta(2(2h-1))^{1/2}}+\frac{\sqrt{2}Gh^2+(1-\xi)D\mu}{\delta\mu(2(2h-1))^{1/4}}\Big)\sum_{t=1}^T\varphi_t  \\
& \quad + \xi G T \sqrt{h} D,\label{eqn:initial step regret}
\end{align}
with $V_{T}\triangleq\E\big[\sum_{t=1}^{T-1}\norm{x_{t}^{\star}-x_{t+1}^{\star}}\big]$, where $\E[\cdot]$ denotes the expectation with respect to the randomness in Algorithm~\ref{alg:two point feedback} and $\bm{x}^\star$ is the global optimal solution defined in \eqref{eqn:offline obj}.
\end{theorem}

Recalling the definition of Algorithm~\ref{alg:two point feedback}, the initialization phase of Algorithm~\ref{alg:two point feedback} can be viewed as being applied to the cost functions $f_1(\cdot),\dots,f_T(\cdot)$ (without using any prediction) and returning $x_1^0,\dots,x_T^0$. To be more specific, the initialization phase of Algorithm~\ref{alg:two point feedback} solves the problem of OCO with memory, where the cost functions are given by $f_t(\cdot),\dots,f_T(\cdot)$ and one can check that by line~7 of Algorithm~\ref{alg:two point feedback}, the decision points $x_1^0,\dots,x_T^0$ are obtained as
\begin{equation*}
x_{t+1}^0=\Pi_{(1-\xi)\X}(x_t^0-\eta_tg_t^0),
\end{equation*}
where $g_t^0$ is given by line~6 of Algorithm~\ref{alg:two point feedback}. As we argued Section~\ref{sec:problem formulation}, our OCO with memory framework corresponds to the bandit information setting, i.e., $g_t^0$ is computed based on only the (noisy) cost function values. The resulting dynamic regret of the initialization step applied to the problem of bandit OCO with memory is then given by \eqref{eqn:initial step regret} in Theorem~\ref{thm:initial regret bound}. The OCO with memory framework has been studied in \cite{arora2012online,anava2015online,zhao2022non}. We extend the existing algorithm design and regret analysis to the bandit OCO with memory framework, where the function value oracle can potentially be noisy. In particular, setting $\varphi_t=0$ for all $t\ge0$ (i.e., there is no noise in the oracle of the function values of $f_t(\cdot)$), and setting  $\delta=1/\sqrt{T}$ and $\xi=\frac{1}{r\sqrt{T}}$ according to the discussions in Remark~\ref{remark:feasible small set}, one can show that the regret bound in \eqref{eqn:initial step regret} reduces to $\tilde{\CO}(V_{T}\sqrt{T})$, where $\tilde{\CO}(\cdot)$ hides logarithmic factors in $T$.  Moreover, one can show that the $\tilde{\CO}(\sqrt{T})$ scaling of the regret still holds when $\varphi_t\le1/T$. Also note that the $V_{T}$ factor in the regret bound in \eqref{eqn:initial step regret} can potentially be refined to $\sqrt{V_{T}}$ using the techniques from \cite{jadbabaie2015online,zhao2022non}. Specifically, supposing $V_{T}$ (or an upper bound on $V_{T}$) is known, then setting the step size $\eta=\frac{1}{t\mu\sqrt{V_T}}$ in Algorithm~\ref{alg:two point feedback} yields a regret bound of $\tilde{\CO}(\sqrt{V_{T} T})$, which matches the regret bound for the algorithm proposed in \cite{zhao2022non} for full information OCO with memory (i.e., the cost functions and gradients are fully accessible) up to logrithmic factors in $T$.

{\bf Zeroth-Order Update Phase.} Next, we upper bound the contribution of the zeroth-order method used in Algorithm~\ref{alg:two point feedback} to $R(\A)$. As we argued before in Section~\ref{sec:algorithm design}, the iterates $\bm{x}^j=(x_1^j,\dots,x_T^j)$ for $j=0,\dots,K-1$ in Algorithm~\ref{alg:two point feedback} are obtained from the update rule $\bm{x}^{j+1}=\mathbf{\Pi}_{(1-\xi^\prime)\X}(\bm{x}^j-\alpha\bm{g}^j)$ starting from the initial point $\bm{x}^0$ (given by the initialization phase of Algorithm~\ref{alg:two point feedback}), where $\bm{g}^j=(\bm{g}^j_1,\dots,\bm{g}^j_T)\in\R^{dT}$ with $\bm{g}^j_s=\sum_{k=s}^{s+h-1}g_k^j$ for all $s\in[T]$ and $\bm{g}^j$ can be viewed as an estimate of the true gradient $\nabla C_T(\bm{x}^j)$ based on the function values of $l_t(\cdot)$ evaluated at certain points perturbed by the truncated Gaussian smoothing (see our arguments in Section~\ref{sec:algorithm design}).

The following result characterizes the total dynamic regret of Algorithm~\ref{alg:two point feedback} by combining the initialization regret and the convergence properties of the zeroth-order update. The detailed proof, including the decomposition of regret and the convergence analysis, is provided in Appendix~\ref{app:proofs of main results}.

\begin{theorem}\label{thm:overall convergence}
Consider the same setting as Theorem~\ref{thm:initial regret bound} and further let $\alpha=1/\beta^{\prime}$ in Algorithm~\ref{alg:two point feedback}, where $\beta^{\prime}=\beta h$ with $\beta$ given by Assumption~\ref{ass:objective functions}. Then, the overall regret $R(\A)$ of Algorithm~\ref{alg:two point feedback} defined in \eqref{eqn:def of regret} satisfies:
\begin{align}\label{eqn:overall convergence of alg}
\E[R(\A)] \le \Big(\frac{1}{1+\gamma}\Big)^K &\E[C_T(\bm{x}^0)-C_T(\bm{x}^\star)] + \frac{\varepsilon}{\gamma} + \xi^\prime G T \sqrt{h} D,
\end{align}
where $\E[C_T(\bm{x}^0)-C_T(\bm{x}^\star)]$ is the initialization regret given in Theorem~\ref{thm:initial regret bound}, $\gamma\triangleq\mu/(\beta^{\prime}-\mu)$\footnote{For any $\beta$-smooth and $\mu$-strongly convex function $f(\cdot)$, it is known that $\beta>\mu$ holds (e.g., \cite{nesterov2013introductory}).}, and
\begin{align*}
     \varepsilon\triangleq\frac{(1-\xi^\prime)D\delta^{\prime}\beta h\sqrt{T}}{2\sqrt{2}(2(2h-1))^{3/4}}+\sqrt{h\beta GT\delta^{\prime}}(1-\xi^\prime)D(Td+3)^{3/4}+\frac{(1-\xi^\prime)Dh\sum_{t=1}^T\varphi_t}{\delta^{\prime}(2(2h-1))^{1/4}}.
\end{align*}
\end{theorem}

We see from Theorem~\ref{thm:overall convergence} that for the zeroth-order update phase (lines~8-20) to be effective, the number of iterations $K$ must be at least $1$. Since $K = \lfloor W / (h-1) \rfloor$, a minimum requirement is that $W \ge h-1$. Note that we do not treat $W$ to be an input parameter to Algorithm~\ref{alg:two point feedback}; rather, the value of $W$ depends on the specific problem setup. In practice, $W$ captures the ability to obtain predictions using the simulator $l_t(\cdot)$.
Specifically, Theorem~\ref{thm:overall convergence} shows that the dynamic regret decays exponentially with $K$, and thus with $W$, which suggests a larger $W$ is always beneficial for achieving lower regret. However, in many real-world applications, obtaining reliable predictions becomes increasingly difficult as the prediction horizon grows, since the quality and availability of future information may degrade significantly for a longer prediction window. 

Moreover, the convergence of the first term on the right-hand side of \eqref{eqn:overall convergence of alg} achieves the same linear convergence rate (i.e., $(1-\mu/\beta^{\prime})^K$) as the standard first-order method based on gradient descent for smooth and strongly convex cost functions \cite{nesterov2013introductory}. Supposing $\varphi_t=0$ for all $t\ge0$, the second term on the right-hand side of \eqref{eqn:overall convergence of alg} can be made arbitrarily small by setting $\delta^{\prime}$. A careful inspection of the proof of Theorem~\ref{thm:overall convergence} reveals that the result in Theorem~\ref{thm:overall convergence} still holds if we remove Assumption~\ref{ass:objective functions}(c). The reason is that in the proof, we only need an upper bound on $\nabla C_T(\bm{x})$ for all $\bm{x}\in\X^{dT}$. Since $C_T(\cdot)$ is $\beta\sqrt{2h}$-smooth from Lemma~\ref{lemma:properties of hat f_t^c}, i.e., $\nabla C_T(\cdot)$ is continuous, and $\X^{dT}$ is bounded by Assumption~\ref{ass:feasible set}, such an upper bound on $\nabla C_T(\cdot)$ always exists. Based on the above arguments, under the smoothness and strong convexity assumption and setting $\varphi_t=0$ for all $t\ge0$, our convergence rate given in Theorem~\ref{thm:overall convergence} for the zeroth-order subroutine in Algorithm~\ref{alg:two point feedback} improves that of the zeroth-order method provided in \cite{nesterov2017random}.\footnote{Since \cite{nesterov2017random} considers the setting of $\varphi_t=0$, we also $\varphi_t=0$ for comparison.} Specifically, the zeroth-order method in \cite{nesterov2017random} also has the form $\bm{x}^{j+1}=\bm{x}^j-\alpha\bm{g}^j$ as we described above, where $\bm{g}^j$ is obtained based on function values at certain points perturbed by a Gaussian smoothing, i.e., $u_s^{j+1}$ in line~12 of Algorithm~\ref{alg:two point feedback} is drawn from a Gaussian distribution $\N(0,I_d)$. Setting the step size to satisfy $\alpha=\frac{1}{4(Td+4)\beta^\prime}$, where $Td$ is the dimension of the variable $\bm{x}^j$, \cite{nesterov2017random} provides a convergence
\begin{align}
\E\big[C_T(\bm{x}^K)-C_T(\bm{x}^{\star})\big]\le\Big(1-\frac{\mu}{8(Td+4)\beta^\prime})\Big)^K\big(C_T(\bm{x}^0)
-C_T(\bm{x}^{\star})\big)+\CO(\frac{\delta^{\prime}T^2d^2\beta^{\prime}}{\mu}).\label{eqn:convergence of gaussian smoothing} 
\end{align}
Comparing \eqref{eqn:convergence of gaussian smoothing} to \eqref{eqn:overall convergence of alg} and noting that $\gamma=\mu/(\beta^{\prime}-\mu)$, we see that Algorithm~\ref{alg:two point feedback} achieves a faster convergence around a neighborhood of the optimal solution $C_T(\bm{x}^{\star})$ and the convergence rate of Algorithm~\ref{alg:two point feedback} does not degrade as the horizon length $T$ increases.

Furthermore, the zeroth-order method based on Gaussian smoothing considered in \cite{nesterov2017random} provides no guarantee that its perturbed query points (e.g., $\bar{x}=x+\delta u$ with $u\sim\N(0,I_d)$) remain within the feasible set $\mathcal{X}$, due to the fact that the norm of the Gaussian random vector $u$ may not be bounded. In contrast, as we argued in Remark~\ref{remark:feasible small set}, 
we can guarantee that the perturbed query points (e.g., $\bar{x}=x+\delta u$ with $u\sim\tilde{N}(0,I_d)$) remain within the feasible set $\mathcal{X}$, by projecting $x$ onto smaller set $(1-\xi)\X$ and leverging the fact that the truncated Gaussian random vector $u$ always has bounded norm. 

{\bf Overall Bound on the Regret $R(\A)$ of Algorithm~\ref{alg:two point feedback}.} Combining the result in Theorems \ref{thm:initial regret bound}-\ref{thm:overall convergence} gives the desired dynamic regret of Algorithm \ref{alg:two point feedback} for solving OCO with memory and limited predictions. The total regret consists of an exponentially decaying term of the initialization regret with respect to $K = \lfloor W/(h-1) \rfloor$, plus the error terms $\varepsilon/\gamma$ and $\xi^\prime G T \sqrt{h} D$. Supposing $\varphi_t=0$ for all $t\ge0$, the two error terms can be made arbitrarily small by setting $\delta^{\prime}$, and the overall regret Algorithm~\ref{alg:two point feedback} reduces to the exponential decay term of the initialization regret, which matches that of the algorithm proposed in \cite{li2020online} for the full information setting with one-step memory (i.e., $h=1$).

\section{Numerical Experiments}\label{sec:simulations}
This section provides numerical results to complement our theoretical analysis and demonstrate the performance of our proposed algorithm.
\subsection{Experimental Setup}\label{subsec:Experimental Setup}
We consider an unconstrained quadratic programming problem, a standard benchmark in convex optimization, where the total cost is given by $C(x)=\sum_{t=1}^T f_t( x_{t-h+1:t})$. The cost at each time step is $f_t(x_{t-h+1:t})=\frac{1}{2}x_{t-h+1:t}^\top A_t x_{t-h+1:t}+B_t^\top x_{t-h+1:t}$. For all experiments, we consider a scalar decision variable $x_t\in\R$ (i.e., $d=1$) and a memory of $h=2$. The matrices $A_t\in\R^{2\times 2}$ and vectors $ B_t\in\R^2$ are randomly generated at each run, with $A_t$ being constructed as positive definite to ensure that the strong convexity and smoothness properties from Assumption~\ref{ass:feasible set} hold. The algorithm parameters are set according to our theoretical analysis: $x_0=0.5,\eta_t=0.2/t,\alpha=0.05,\delta=0.2,\delta^{\prime}=0.0001$. We assume a perfect simulator, i.e., $l_t(x)=f_t(x)$ for all $t$ ($\phi_t=0$). We evaluate the performance of our algorithm under three different smoothing distributions used for the two-point feedback gradient estimation:
\begin{itemize}
\item[(a)] 
\textbf{Gaussian distribution:} The random direction $u$ is drawn from a standard normal distribution, $u\sim\N(0,1)$.
\item[(b)] \textbf{Truncated Gaussian distribution:} The random direction $u$ is drawn from a normal distribution truncated to the interval $[-2,2]$, which reduces the variance of the gradient estimate.
\item[(c)] \textbf{Bernoulli distribution:} The random direction $u$ is drawn from a Bernoulli distribution, taking values in $\{-1,1\}$ with equal probability.
\end{itemize}
All results are generated by averaging over 100 independent experimental runs. The plots show the median regret, with the shaded areas representing the interquartile range (from the 25th to the 75th percentile).
\subsection{Validation of Initialization Phase Regret (Theorem~\ref{thm:initial regret bound})}
\begin{figure}[ht]
\begin{subfigure}{0.49\textwidth}
  \centering
  \includegraphics[width=0.99\linewidth]{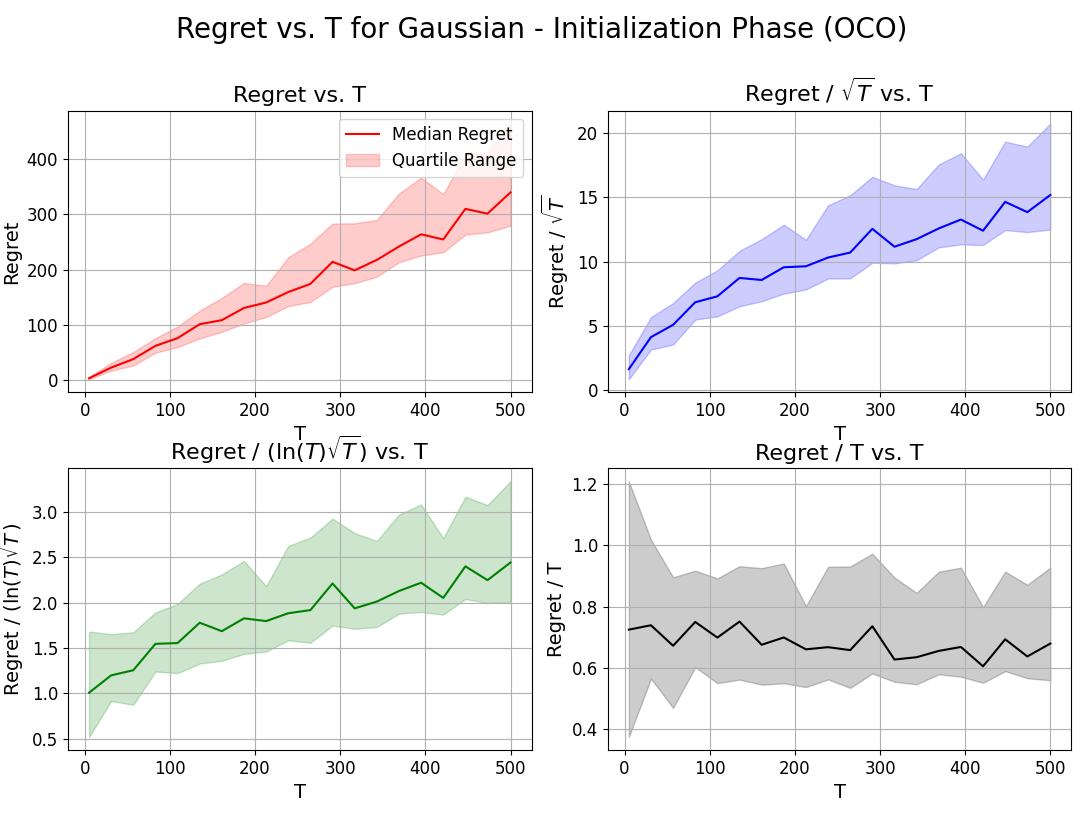}
  \subcaption{Gaussian}
\end{subfigure}
\hfill
\begin{subfigure}{0.49\textwidth}
  \centering
  \includegraphics[width=0.99\linewidth]{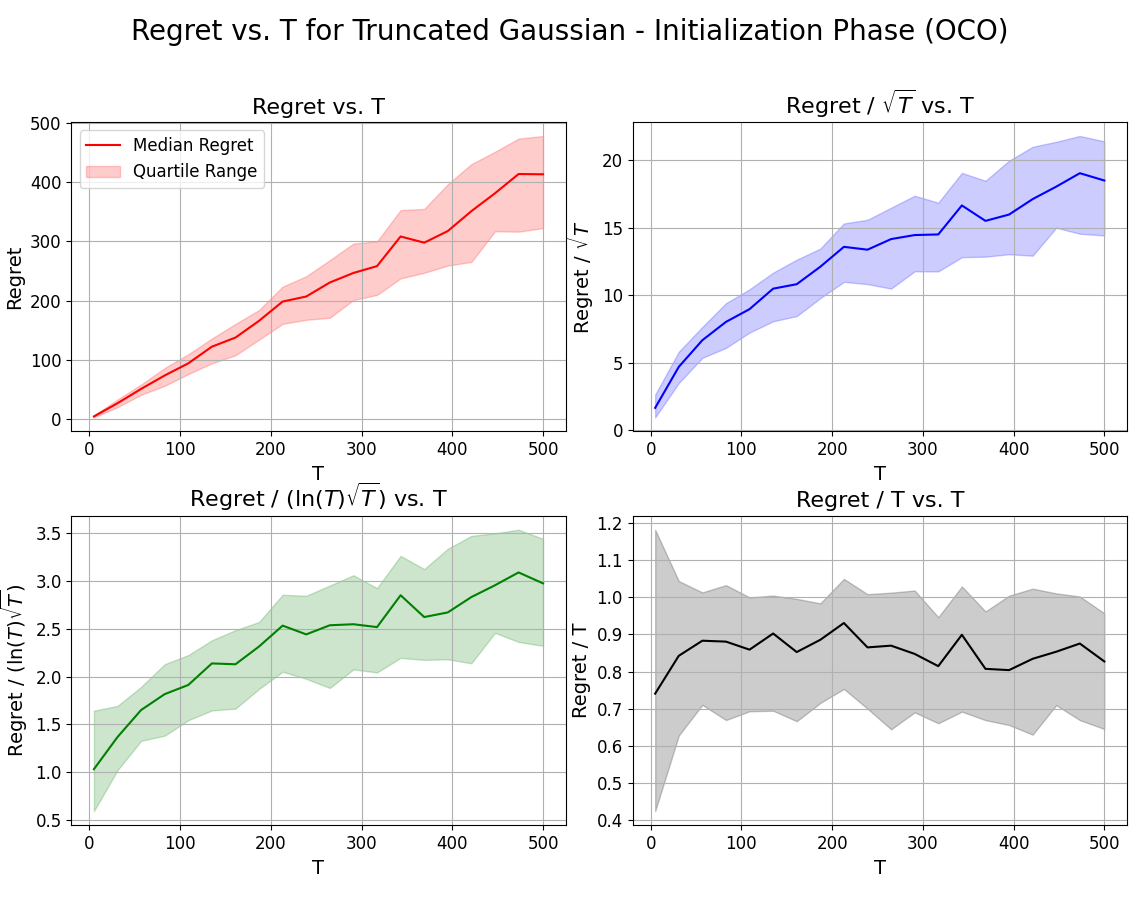}
  \subcaption{Truncated Gaussian}
\end{subfigure}
\hfill
\newline
\centering
\begin{subfigure}{0.49\textwidth}
  \centering
  \includegraphics[width=0.99\linewidth]{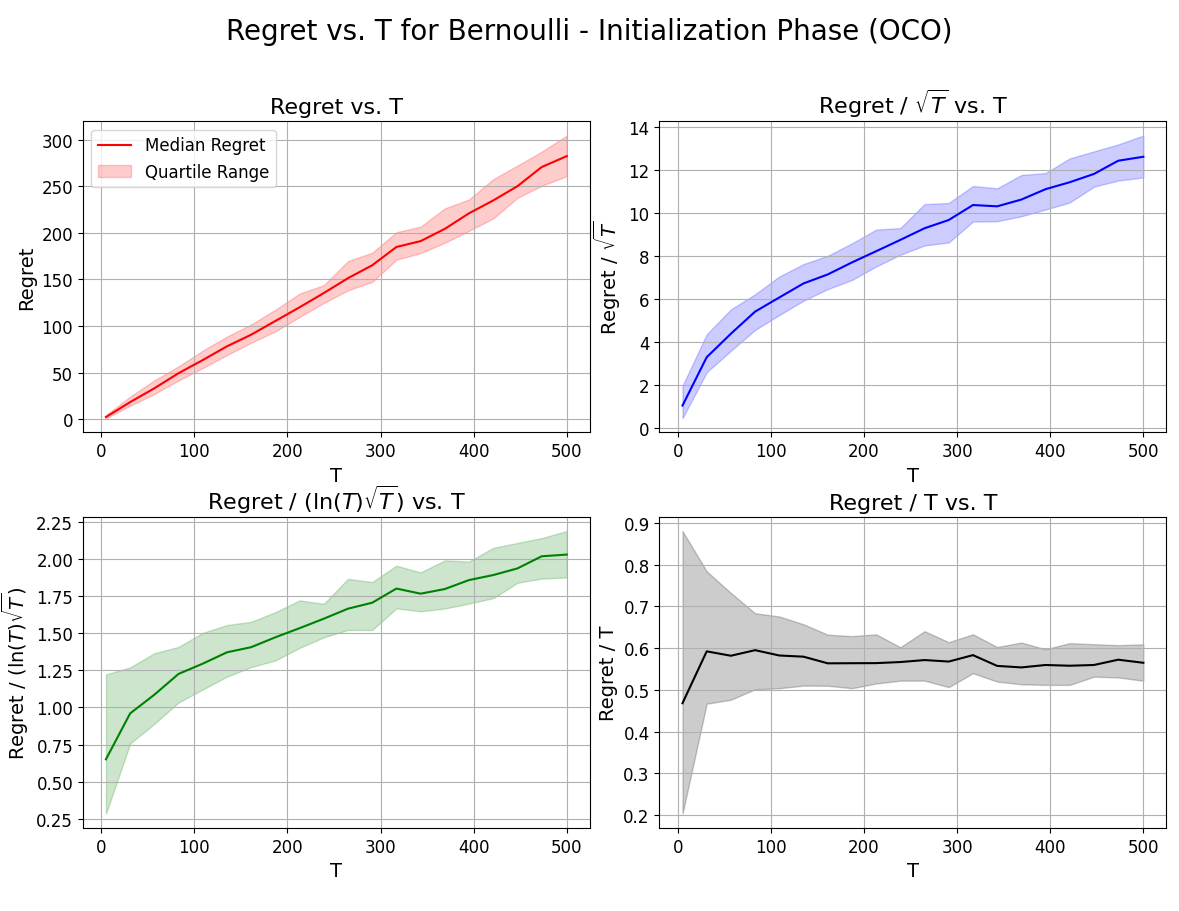}
  \subcaption{Bernoulli}
\end{subfigure}
\hfill
\caption{
}
\label{fig:init_regret}
\end{figure}

First, we empirically validate the performance of the initialization phase of our algorithm (lines 3-8 in Algorithm~\ref{alg:two point feedback}), which is an Online Convex Optimization (OCO) procedure with memory. Theorem~\ref{thm:initial regret bound} establishes a dynamic regret bound of $\O(TlnT)$ for this phase. To verify this, we run the initialization phase over a long time horizon, varying $T$ from 5 to 500, and compute the regret of the generated sequence $x_0$. Figure~\ref{fig:init_regret} displays the normalized regret for all three smoothing distributions. The subplots show that the regret grows sub-linearly with $T$, as evidenced by the decreasing curves in the $Reg/T$ plots. Furthermore, the $Reg/T$ plots show a slowly growing trend, consistent with the theoretical $\ln T$ factor. These results provide strong empirical support for the regret bound established in Theorem~\ref{thm:initial regret bound}. Among the smoothing methods, the lower-variance distributions (Truncated Gaussian and Bernoulli) tend to exhibit slightly better and more stable performance than the standard Gaussian.
\graphicspath{ {./images/} }

\subsection{Regret Decay with Prediction Window (Theorem~\ref{thm:overall convergence})}
\begin{figure}[htbp]\centering
\includegraphics[width=0.49\textwidth]{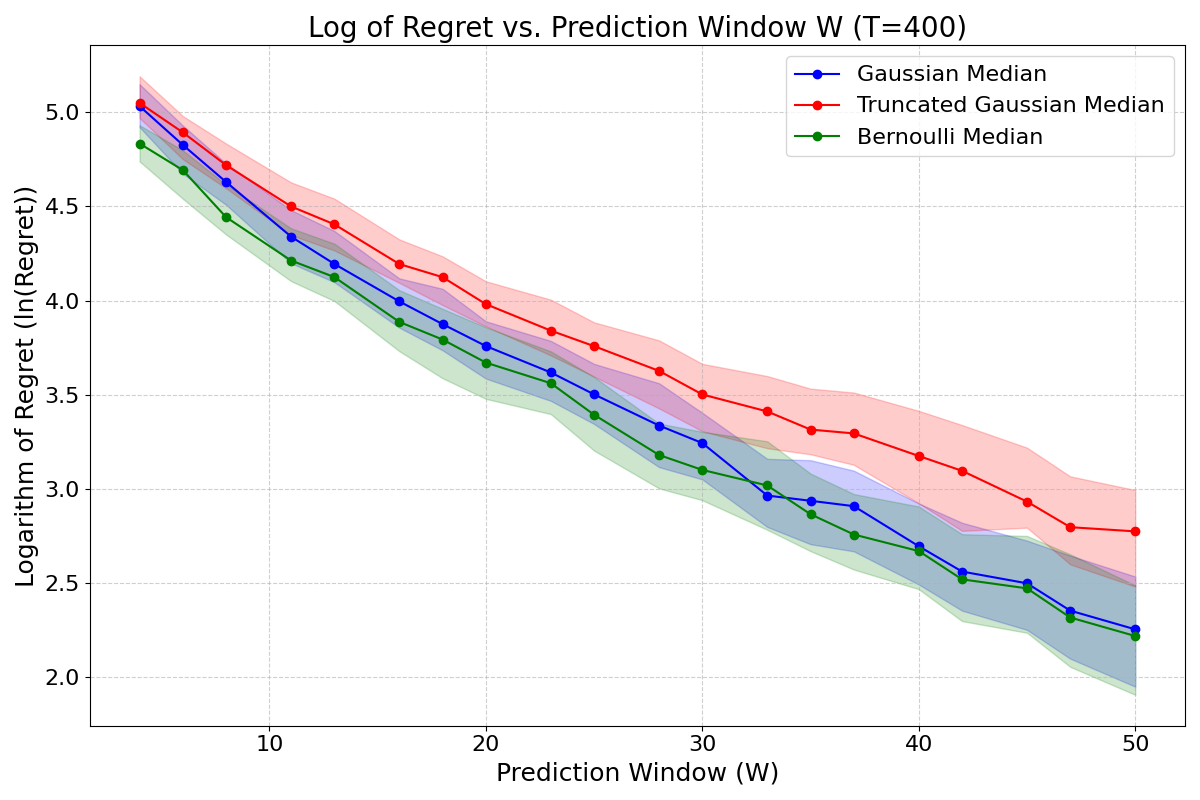}
\caption{Log of regret for the full algorithm versus the prediction window size $W$, with $T=400$. The linear downward trend for all smoothing methods empirically confirms the theoretical exponential decay of regret, a key result derived from Theorem~\ref{thm:overall convergence}.}
\label{fig:regret_vs_W}
\end{figure}

Next, we investigate the core claim of our paper: the total regret of the full algorithm decays exponentially with the length of the prediction window, $W$. This behavior is a direct consequence of the linear convergence of the zeroth-order phase, as analyzed in Theorem~\ref{thm:overall convergence}. For this experiment, we fix the time horizon at $T=400$ and vary the prediction window $W$ from 4 to 50. The regret is calculated based on the final output sequence $x_K$, where $K=\lfloor W/(h-1)\rfloor$ is the number of zeroth-order iterations. Figure~\ref{fig:regret_vs_W} plots the natural logarithm of the regret, $\ln(Regret)$, against $W$. As depicted, all three smoothing methods show a clear, approximately linear decrease on this semi-log plot. This linear trend provides strong empirical evidence of the exponential decay of regret with respect to the prediction window size. This result validates the convergence analysis of Theorem~\ref{thm:overall convergence}, demonstrating that as more future information is made available (larger $W$), the algorithm performs more optimization iterations and rapidly converges to a better solution, drastically reducing the overall regret. Finally, the performances of using the Bernoulli random variable and the Gaussian random variable for smoothing are comparable to the case of using the truncated Gaussian random variable; however, there are no theoretical guarantees of the regret for Bernoulli smoothing or Gaussian smoothing, which are left to future work.

\subsection{Convergence Rate Comparison with \cite{nesterov2017random} } 
\begin{figure}[htbp]
\centering
\includegraphics[width=0.45\textwidth]{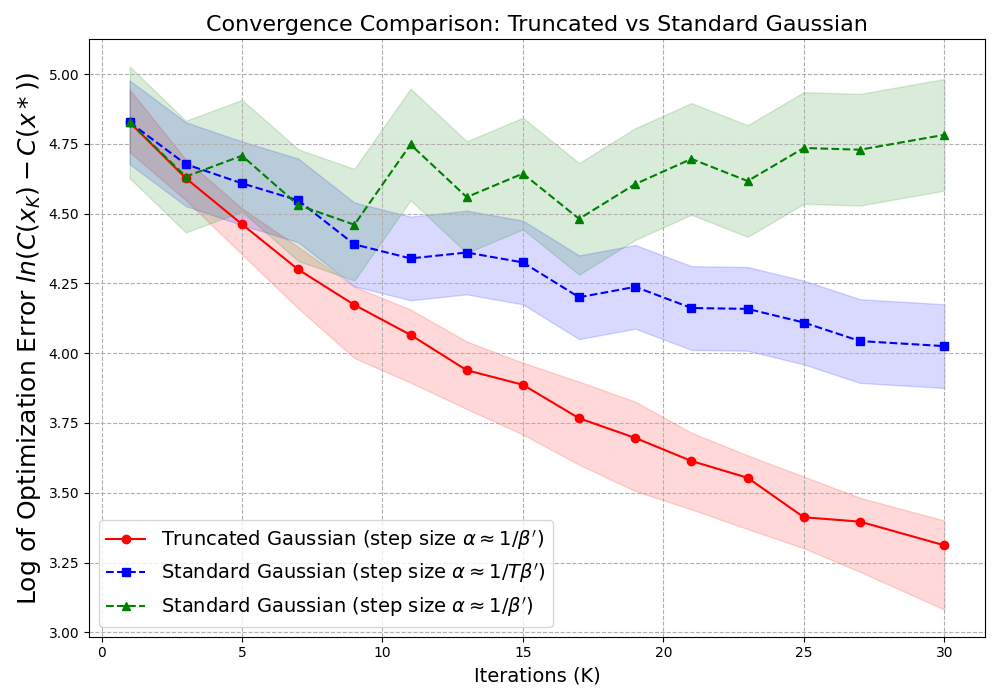}
\caption{Convergence Rate Comparison (T=200)}\label{fig:convergence_rate}
\end{figure}

Finally, to highlight the advantages of our zeroth-order method based on truncated Gaussian smoothing, we directly compare its convergence rate against the zeroth-order method proposed by \cite{nesterov2017random} based on Gaussian smoothing. In this experiment, we fix the horizon to be $T=200$. Both algorithms start from the same initial sequence $x_0$, generated by our OCO initialization phase. Figure~\ref{fig:convergence_rate} aggregates the performance curves under different step-size settings into a single plot for direct comparison.
As we argued in Section~\ref{sec:regret analysis}, our algorithm allows for a large step size $\alpha=1/\beta^{\prime}$ compared to the step size $\alpha=1/(T\beta^{\prime})$ (here $\beta^\prime=20$ according to the setup in Section~\ref{subsec:Experimental Setup}) that can be chosen for the algorithm in \cite{nesterov2017random}, which aligns with the results in Figure~\ref{fig:convergence_rate}. 
Moreover, when choosing the step size to be $\alpha=1/\beta^{\prime}$ for the algorithm in \cite{nesterov2017random}, we see from Figure~\ref{fig:convergence_rate} that the algorithm does not converge. This divergence occurs because the method in \cite{nesterov2017random} theoretically requires a step size scaling as $\mathcal{O}(1/T)$ to ensure stability, whereas our method robustly supports a much larger step size (scaling as $\mathcal{O}(1)$). 
These results empirically validate that our approach not only achieves faster convergence but also offers better stability under a more aggressive step size.

\section{Conclusion}\label{sec:conclusion}
We studied the problem of OCO with memory and predictions over a horizon $T$. We proposed an algorithm to solve this problem and showed that the dynamic regret of the algorithm decays exponentially with the prediction window length. The first subroutine of the algorithm can be used to solve general online convex optimization with memory and bandit feedback with $\sqrt{T}$-dynamic regret with respect to $T$. The second subroutine of the algorithm is a zeroth-order method that can be used to solve general convex optimization problems with a linear convergence rate that improves the rate of zeroth-order algorithms and matches the rate of first-order algorithms for convex optimization problems. Our algorithm design and analysis rely on truncated Gaussian smoothing when querying the decision points for obtaining the predictions. Future work includes extending the algorithm design and analysis to more general (i.e., nonconvex or nonsmooth) cost functions and exploring other possible smoothing methods for zeroth-order algorithms.

\bibliography{main}
\bibliographystyle{unsrt}

\appendix
\counterwithin{lemma}{section}
\counterwithin{theorem}{section}
\counterwithin{proposition}{section}
\counterwithin{corollary}{section}
\counterwithin{definition}{section}
\counterwithin{equation}{section}
\onecolumn

\section{Proofs of Preliminary Results}\label{app:preliminary proofs}
\subsection{Proof of Lemma~\ref{lemma:properties of hat f_t^c}}
We first prove the properties of $\hat{f}_t^c(\cdot)$. Note that the $\mu$-strongly convexity of $f^c_t(\cdot)$ follows directly from that of $f_t(\cdot)$. Thus, for any $x,y\in\R^d$, 
\begin{equation*}
f_t^c(x+\delta v)-f_t^c(y+\delta v)\le \nabla f_t^c(x+\delta v)^{\top}(x-y)-\frac{\mu}{2}\norm{x-y}^2.
\end{equation*}
Taking the expectation with respect to $v\sim\widetilde{\N}(0,I_d)$ on both sides of the above inequality, we get
\begin{align}\nonumber
\hat{f}_t^c(x)-\hat{f}_t^c(y)&\le\E_{v}\big[\nabla f_t^c(x+\delta v)\big]^{\top}(x-y)-\frac{\mu}{2}\norm{x-y}^2\\
&\le\nabla\hat{f}_t^c(x)^{\top}(x-y)-\frac{\mu}{2}\norm{x-y}^2,\label{eqn:strong convex of hat f_t^c}
\end{align}
where we can swap the order of $\E_{v}$ and $\nabla$ since $\nabla f_t^c(\cdot)$ is continuous in $x$. It follows from \eqref{eqn:strong convex of hat f_t^c} that $\hat{f}_t^c(\cdot)$ is $\mu$-strongly convex. Similarly, the $\beta$-smoothness of $f_t^c(\cdot)$ follows from that of $f_t(\cdot)$. For any $x,y\in\R^d$, we have
\begin{equation*}
\norm{\nabla f_t^c(x+\delta v)-\nabla f_t^c(y+\delta v)}\le \beta\norm{x-y}.
\end{equation*}
Again, taking the expectation $\E_v[\cdot]$ on both sides of the above inequality yields
\begin{align}\nonumber
&\E_{v}\big[\norm{\nabla f_t^c(x+\delta v)-\nabla f_t^c(y+\delta v)}\big]\le\beta\norm{x-y},\\\nonumber
\overset{(a)}\Rightarrow &\norm{\E_{v}[\nabla f_t^c(x+\delta v)]-\E_v[\nabla f_t^c(y+\delta v)]}\le\beta\norm{x-y},\\\nonumber
\Rightarrow&\norm{\nabla\hat{f}_t^c(x)-\nabla\hat{f}_t^c(y)}\le\beta\norm{x-y},
\end{align}
where (a) follows from Jensen's inequality (i.e., $\E[\norm{a}]\ge\norm{\E[a]}$ for any $a\in\R^d$), which shows that $\hat{f}_t^c(\cdot)$ is $\beta$-smooth.

Next, we prove the properties of $C_T(\cdot)$. Recall that $C_T(\bm{x})=\sum_{t=1}^Tf_t(x_{t-h+1:t})$, where $\bm{x}=(x_1,\dots,x_T)\in\R^{dT}$ and $x_t=x_0$ for all $t\le0$. Considering any $\bm{x}=(x_1,\dots,x_T)\in\R^{dT}$ and any $\bm{y}=(y_1,\dots,y_T)\in\R^{dT}$, we have 
\begin{equation}
\norm{\bm{x}-\bm{y}}^2\le\sum_{t=1}^T\norm{x_{t-h+1:t}-y_{t-h+1:t}}^2.\label{eqn:upper bound on x-y}
\end{equation}
Moreover, we have
\begin{align}
&\sum_{t=1}^T\norm{x_{t-h+1:t}-y_{t-h+1:t}}^2\le h\norm{\bm{x}-\bm{y}}^2,\label{eqn:lower bound on (x-y)^2} \\ 
\Rightarrow &\sqrt{\sum_{t=1}^T\norm{x_{t-h+1:t}-y_{t-h+1:t}}^2}\le\sqrt{h}\norm{\bm{x}-\bm{y}},\nonumber \\
\overset{(a)}\Rightarrow &\sum_{t=1}^T\norm{x_{t-h+1:t}-y_{t-h+1:t}}\le\sqrt{Th}\norm{\bm{x}-\bm{y}},\label{eqn:lower bound on x-y}
\end{align}
where (a) uses the fact that $\sqrt{\sum_{t=1}^Ta_t^2}\ge\frac{\sum_{t=1}^T a_t}{\sqrt{T}}$ for any $a_t\in\R$. By the $\mu$-strongly convexity of $f_t(\cdot)$, for any $t\in[T]$,
\begin{align}
f_t(x_{t-h+1:t})-f_t(y_{t-h+1:t})\le \nabla f_t(x_{t-h+1:t})^{\top}(x_{t-h+1:t}-y_{t-h+1:t})-\frac{\mu}{2}\norm{x_{t-h+1:t}-y_{t-h+1:t}}^2.\label{eqn:stongly convexity of f_t}
\end{align}
To relate $\nabla f_t(x_{t-h+1:t})$ to $\nabla C_T(\bm{x})$, we first augment the gradient vector $\nabla f_t(x_{t-h+1:t})\in\R^{dh}$ to a higher dimensional vector $\tilde{\nabla}f_t(x_{t-h+1:t})\in\R^{dT}$ by adding zero entries. By the definition of $C_T(\cdot)$, we see that $\nabla C_T(\bm{x})=\sum_{t=1}^T\tilde{\nabla} f_t(x_{t-h+1:t})$. 
Summing \eqref{eqn:stongly convexity of f_t} over $t\in[T]$ now gives us
\begin{align*}
&C_T(\bm{x})-C_T(\bm{y})\le\sum_{t=1}^T\tilde{\nabla}f_t(x_{t-h+1:t})^{\top}(\bm{x}-\bm{y})-\frac{\mu}{2}\sum_{t=1}^T\norm{x_{t-h+1:y}-y_{t-h+1:t}}^2,\\
\overset{\eqref{eqn:upper bound on x-y}}\Rightarrow&C_T(\bm{x})-C_T(\bm{y})\le\nabla C_T(\bm{x})^{\top}(\bm{x}-\bm{y})-\frac{\mu}{2}\norm{\bm{x}-\bm{y}}^2,
\end{align*}
which proves that $C_T(\cdot)$ is $\mu$-strongly convex. To show the smoothness of $C_T(\cdot)$, we again start from $f_t(\cdot)$ and its $\beta$-smoothness gives us\footnote{It is well known that a $\beta$-smooth function $f:\R^{d}\to\R$ as per definition~\ref{def:beta-smooth} satisfies that $ f(y)-f(x)\le \left \langle \nabla f(x),y-x \right \rangle +\frac{\beta}{2}\norm{y-x}^2,\ \forall x,y\in\R^d$ \cite{bubeck2015convex}.} 
\begin{equation}
f_t(y_{t-h+1:t})-f_t(x_{t-h+1:t})\le \nabla f_t(x_{t-h+1:t})^{\top}(y_{t-h+1:t}-x_{t-h+1:t})+\frac{\beta}{2}\norm{x_{t-h+1:t}-y_{t-h+1:t}}^2.\label{eqn:smoothness of f_t}
\end{equation}
Summing \eqref{eqn:smoothness of f_t} over $t\in[T]$ gives us 
\begin{align}\nonumber
&C_T(\bm{y})-C_T(\bm{x})\le \sum_{t=1}^T\tilde{\nabla}f_t(x_{t-h+1:t})^{\top}(\bm{y}-\bm{x})+\frac{\beta}{2}\sum_{t=1}^T\norm{x_{t-h+1:y}-y_{t-h+1:t}}^2,\\ \nonumber
\overset{\eqref{eqn:lower bound on (x-y)^2}}\Rightarrow&C_T(\bm{y})-C_T(\bm{x})\le\nabla C_T(\bm{x})^{\top}(\bm{y}-\bm{x})+\frac{\beta h}{2}\norm{\bm{x}-\bm{y}}^2,
\end{align}
which proves that $C_T(\cdot)$ is $(\beta h)$-smooth. To show that that $C_T(\cdot)$ is $G\sqrt{Th}$-Lipschitz, we have
\begin{align}
|C_T(\bm{x})-C_T(\bm{y})|\le \sum_{t=1}^T|f_t(x_{t-h+1:t})-f_t(y_{t-h+1:t})|\overset{\eqref{eqn:lower bound on x-y}}\le G\sqrt{Th}\norm{\bm{x}-\bm{y}}.
\end{align}

Finally, recalling the definition of $\hat{C}_T(\cdot)$, the $\mu$-strong convexity and $(\beta h)$-smoothness of $\hat{C}_T(\cdot)$ follows from similar arguments to those above for $\hat{f}_t^c(\cdot)$.$\hfill\blacksquare$

\subsection{Proof of Lemma~\ref{lemma:zeroth-order gradient}}
{\it Proof of part (a)} Eq.~\eqref{eqn:two point estimator for hat f_t^c} follows directly from \citep[Eq.~(26)]{nesterov2017random}. Since $f_t^c(\cdot)$ is $\beta$-smooth by Assumption~\ref{ass:objective functions}(a), Eq.~\eqref{eqn:f_t^c and hat f_t^c distance} follows from \citep[Theorem~1]{nesterov2017random}.

\noindent{\it Proof of part (b)} Recalling the definition of $\hat{C}_T(\cdot)$, we have 
\begin{equation}
\hat{C}_T(\bm{x})=\sum_{t=1}^T\E_{\bm{v}}[f_t(\bar{x}_{t-h+1:t})],
\end{equation}
where we set $\bar{x}_t=x_0$ for all $t\le0$. It follows that for any $s\in[T]$, 
\begin{align}\nonumber
\frac{\partial\hat{C}_T(\bm{x})}{\partial x_s}&=\frac{\partial\E_{\bm{v}}\big[\sum_{t=1}^Tf_t(\bar{x}_{t-h+1:t})\big]}{\partial x_s}\\\nonumber
&\overset{(a)}=\E_{\bm{v}}\Big[\sum_{t=1}^T\frac{\partial f_t(\bar{x}_{t-h+1:t})}{\partial x_s}\Big]\\\nonumber
&=\E_{\bm{v}}\Big[\sum_{k=s}^{s+h-1}\frac{\partial f_k(\bar{x}_{k-h+1:k})}{\partial x_s}\Big]\\
&\overset{(b)}=\frac{\partial\E_{\bm{v}}\big[\sum_{k=s}^{s+h-1}f_k(\bar{x}_{k-h+1:k})\big]}{\partial x_s},\label{eqn:partial x_s of hat C_T}
\end{align}
where we can swap the order of $\E_{\bm{v}}$ and $\frac{\partial}{\partial x_s}$ in (a) and (b) since $\frac{\partial f_t(\bar{x}_{t-h+1:t})}{\partial x_s}$ is continuous in $x_s$ for all $t\in[T]$. Now, we notice that $\sum_{k=s}^{s+h-1}f_k(x_{k-h+1:k})$  may be viewed as a function of $x_{s-h+1},\dots,x_{s+h-1}$, denoted as $F_s(x_{s-h+1:s+h-1})$. Since $\bm{v}=(v_1,\dots,v_T)\sim\widetilde{\N}(0,I_{dT})$, we have $v_t\sim\widetilde{\N} (0,I_{d})$ for all $t\in[T]$. Letting $\bm{v}_s=v_{s-h+1:s+h-1}\in\R^{2h-1}$ be a column vector and also recalling that $\bar{x}_t=x_t+\delta^{\prime}v_t$ and $\tilde{x}_t=x_t-\delta^{\prime}v_t$, we again have from \citep[Eq.~(26)]{nesterov2017random} that 
\begin{align}\nonumber
\nabla\E_{\bm{v}_s}[F_s(\bar{x}_{s-h+1:s+h-1})]&=\E_{\bm{v}_s}\Big[\frac{F(\bar{x}_{s-h+1:s+h-1})-F(\tilde{x}_{s-h+1:s+h-1})}{2\delta^{\prime}}\bm{v}_s\Big]\\
&=\E_{\bm{v}_s}\Big[\frac{\sum_{k=s}^{s+h-1}\big(f_k(\bar{x}_{k-h+1:k})-f_k(\tilde{x}_{k-h+1:k})}{2\delta^{\prime}}\bm{v}_s\Big].\label{eqn:gradient of exp F_s}
\end{align}
Noting that $\partial\E_{\bm{v}}\big[\sum_{k=s}^{s+h-1}f_k(\bar{x}_{k-h+1:k})\big]/\partial x_s=\partial\E_{\bm{v}_s}[F_s(\bar{x}_{s-h+1:s+h-1})]/\partial x_s$, we see from \eqref{eqn:partial x_s of hat C_T} and \eqref{eqn:gradient of exp F_s} that the partial derivative $\partial \hat{C}_T(\bm{x})/\partial x_s$ can be extracted from the gradient $\nabla\E_{\bm{v}_s}[F_s(\bar{x}_{s-h+1:s+h-1})]$ using \eqref{eqn:two point estimator for hat C_T}. Finally, Eq.~\eqref{eqn:C_T and hat C_T difference} follows from the $\beta h$-smoothness of $C_T(\cdot)$ (by Lemma~\ref{lemma:properties of hat f_t^c}) and \citep[Theorem~1]{nesterov2017random}, and Eq.~\eqref{eqn:C_T and hat C_T gradient difference} follows from the $\beta h$-smoothness of $C_T(\cdot)$ and \citep[Lemma~3]{nesterov2017random}. $\hfill\blacksquare$

\section{Proofs of Main Results}\label{app:proofs of main results}
\subsection{Proof of Theorem~\ref{thm:initial regret bound}}
As outlined in the main text, the regret of the initialization phase, $R(\A_\text{init})$, is defined against the true optimal solution $\bm{x}^{\star} = (x_1^\star, \dots, x_T^\star)$.
To avoid dependencies on the specific optimum of the constrained set $(1-\xi)\mathcal{X}$, we construct a comparison sequence $\hat{\bm{x}}^\star = (\hat{x}_1^\star, \dots, \hat{x}_T^\star)$ by scaling the global optimal solution:
\begin{equation}\label{eqn:def of hat_x_star}
    \hat{x}_t^\star = (1-\xi)x_t^\star, \quad \forall t \in [T].
\end{equation}
Since we assume the feasible set $\mathcal{X}$ contains the origin (or is centered such that this scaling is valid), and $\mathcal{X}$ is convex, it holds that $\hat{x}_t^\star \in (1-\xi)\mathcal{X}$.
We decompose the total regret as follows:
\begin{align*}
\E[R(\A_\text{init})] &= \E\Big[C_T(\bm{x}^{0}) - C_T(\bm{x}^{\star})\Big] \\
&= \underbrace{\E\Big[ C_T(\bm{x}^{0}) - C_T(\hat{\bm{x}}^\star) \Big]}_{\text{Part 1: Regret against scaled comparator}} + \underbrace{\Big( C_T(\hat{\bm{x}}^\star) - C_T(\bm{x}^{\star}) \Big)}_{\text{Part 2: Approximation Error}}.
\end{align*}
We will now bound each part separately.

{\bf Part 1: Bounding the Regret against $\hat{\bm{x}}^\star$.} This part bounds the term $\E\big[ C_T(\bm{x}^{0}) - C_T(\hat{\bm{x}}^\star) \big]$. The proof proceeds in several steps.

{\it step 1: Preparatory Terminology.}
Based on the definition of Algorithm~\ref{alg:two point feedback}, let us further define an auxiliary function 
\begin{equation*}
h_t(x)=\hat{f}_t^c(x)+(\tilde{g}_t-\nabla\hat{f}_t^c(x_t^0))x,
\end{equation*}
where $\tilde{g}_t=\frac{1}{2\delta}\big[f_t^c(\bar{x}_t^0)-f^c_t(\tilde{x}_t^0)\big]u_t^0$. The above definitions imply $\nabla h_t(x_t^0)=\tilde{g}_t$. Moreover, we have from Lemma~\ref{lemma:zeroth-order gradient}(a) that $\E_{u_t}[\tilde{g}_t]=\nabla\hat{f}_t^c(x_t^0)$ and thus $\E_{u_t}[h_t(x)]=\hat{f}_t^c(x)$ (for all $x$ that is independent of $u_t$), which further implies that $\E[h_t(x_t^0)]=\E[\hat{f}_t^c(x_t^0)]$ for all $t$. For our analysis in the sequel, we show that for any $t\in[T]$,
\begin{align}\nonumber
\norm{\nabla h_t(x_t^0)}&=\big\lVert\frac{1}{2\delta}\big[f_t^c(\bar{x}_t^0)-f_t^c(\tilde{x}_t^0)\big]u^0\big\rVert\\\nonumber
&\le\frac{1}{2\delta}|f_t^c(\bar{x}_t^0)-f_t^c(\tilde{x}_t^0)|\norm{u^0}\\\nonumber
&\overset{(a)}\le\frac{1}{2\delta}G\sqrt{h}\norm{\bar{x}_t^0-\tilde{x}_t^0}\norm{u^0}\\
&=\frac{1}{2\delta}G\sqrt{h}\norm{2\delta u^0}\norm{u^0}=G\sqrt{h}\norm{u^0}^2,\label{eqn:upper bound on tilde g_t}
\end{align}
where $(a)$ follows from the $G$-Lipschitz of $f_t(\cdot)$ in  Assumption~\ref{ass:objective functions}(c) and the definition of $f_t^c(\cdot)$. 
Let us denote  $e_t=e_t^m+e_t^c$ with 
\begin{align}
&e_t^m = \frac{1}{2\delta}\Big[\big(l_t(\bar{x}_{t-h+1:t}^0)-l_t(\tilde{x}_{t-h+1:t}^0)\big)-\big(f_t(\bar{x}_{t-h+1:t}^0)-f_t(\tilde{x}_{t-h+1:t}^0)\big)\Big]u^0,\label{eqn:def of e_t^m}\\
&e_t^c = \frac{1}{2\delta}\Big[\big(f_t(\bar{x}_{t-h+1:t}^0)-f_t(\tilde{x}_{t-h+1:t}^0)\big)-\big(f_t^c(\bar{x}_{t}^0)-f_t^c(\tilde{x}_{t}^0)\big)\Big]u^0.\label{eqn:def of e_t^c}
\end{align}

{\it Step 2: Relate Total Cost $C_T(\cdot)$ to Unitary Cost $f_t^c(\cdot)$.} First, we relate the algorithmic regret in terms of $C_T(\cdot)$ to the regret in terms of the $h$-unitary function $f_t^c(x) = f_t(x, \dots, x)$.
\begin{align*}
C_T(\bm{x}^{0}) - C_T(\hat{\bm{x}}^\star) = \sum_{t=1}^T \Big( f_t(x_{t-h+1:t}^0) - f_t(\hat{x}_{t-h+1:t}^\star) \Big).
\end{align*}
We can rewrite this as:
\begin{align*}
\sum_{t=1}^T \Big( f_t^c(x_t^0) - f_t^c(\hat{x}_t^\star) \Big) + \sum_{t=1}^T \Big( f_t(x_{t-h+1:t}^0) - f_t^c(x_t^0) \Big) - \sum_{t=1}^T \Big( f_t(\hat{x}_{t-h+1:t}^\star) - f_t^c(\hat{x}_t^\star) \Big).
\end{align*}
We bound the expectation of the second term using the Lipschitz property (Assumption~\ref{ass:objective functions}(c)):
\begin{align}\nonumber
&\quad\E\Big[\Big|\sum_{t=1}^Tf_t(x_{t-h+1:t}^0)-f_t^c(x_t^0)\Big|\Big]\\\nonumber
&\le G\sum_{t=1}^T\sum_{j=1}^{h-1}\E\big[\norm{x_t^0-x_{t-j}^0}\big] \le G\sum_{t=1}^T\sum_{j=1}^{h-1}\sum_{l=1}^j\E\big[\norm{x_{t-l+1}^0-x_{t-l}^0}\big]\\\nonumber
&\overset{(a)}\le G\sum_{t=1}^T\sum_{j=1}^{h-1}\sum_{l=1}^j\eta_{t-l}\E\big[\norm{\nabla h(x_{t-l}^0)}+\norm{e_{t-l}}\big]\\\nonumber
&\overset{(b)}\le \frac{1}{\sqrt{2(2h-1)}}h^{5/2}G^2\sum_{t=1}^T\eta_t+Gh^2\sum_{t=1}^T\eta_t\norm{e_t}\\\nonumber
&\overset{(c)}\le\frac{h^{5/2}G^2}{\sqrt{2(2h-1)}\mu}(1+\ln T)+\frac{\sqrt{2}Gh^2}{\delta\mu(2(2h-1))^{1/4}}\sum_{t=1}^T\varphi_t+\frac{2G^2h^3}{\mu(2(2h-1))^{1/2}}(1+\ln T)\\
&\le\frac{3G^2h^3}{\mu(2(2h-1))^{1/2}}(1+\ln T)+\frac{\sqrt{2}Gh^2}{\delta\mu(2(2h-1))^{1/4}}\sum_{t=1}^T\varphi_t, \label{eqn:app C8}
\end{align}
where (a) follows from Lemma~\ref{lemma:pythagorean}, (b) from \eqref{eqn:upper bound on tilde g_t} and $\lVert u^0\rVert^2\le\sqrt{1/2(2h-1)}$, and (c) from \eqref{eqn:upper bound on sum eta_t time e_t}.
Then, we bound the third term (which is deterministic) similarly:
\begin{align}\nonumber
\Big|\sum_{t=1}^T f_t(\hat{x}_{t-h+1:t}^\star)-f_t^c(\hat{x}_t^\star)\Big|&\le G\sum_{t=1}^T\sum_{j=1}^{h-1}\sum_{l=1}^j\norm{\hat{x}_{t-l+1}^{\star}-\hat{x}_{t-l}^{\star}}\\\nonumber
&\le Gh^2\sum_{t=1}^T\norm{\hat{x}_t^\star-\hat{x}_{t-1}^{\star}}\\
&\overset{\eqref{eqn:def of hat_x_star}}=(1-\xi)Gh^2\sum_{t=1}^{T-1}\norm{x_t^\star-x_{t+1}^{\star}} = (1-\xi)Gh^2 V_{T},\label{eqn:app C9}
\end{align}
where $V_{T} = \sum_{t=1}^{T-1}\norm{x_t^\star-x_{t+1}^{\star}}$ is the path-length discussed before.

{\it Step 3: Relate Unitary Cost $f_t^c(\cdot)$ to Smoothed Cost $\hat{f}_t^c(\cdot)$.}
Next, we relate $f_t^c(\cdot)$ to its smoothed version $\hat{f}_t^c(\cdot)$ using Lemma~\ref{lemma:zeroth-order gradient}(a), equation \eqref{eqn:f_t^c and hat f_t^c distance}:
\begin{equation*}
|f_t^c(x)-\hat{f}_t^c(x)|\le\frac{\delta^2}{2}\beta d,
\end{equation*}
which implies:
\begin{align*}
f_t^c(x_t^0) - f_t^c(\hat{x}_t^\star) &= \big(f_t^c(x_t^0) - \hat{f}_t^c(x_t^0)\big) + \big(\hat{f}_t^c(x_t^0) - \hat{f}_t^c(\hat{x}_t^\star)\big) + \big(\hat{f}_t^c(\hat{x}_t^\star) - f_t^c(\hat{x}_t^\star)\big) \\
&\le \frac{\delta^2\beta d}{2} + \big(\hat{f}_t^c(x_t^0) - \hat{f}_t^c(\hat{x}_t^\star)\big) + \frac{\delta^2\beta d}{2}.
\end{align*}
Summing over $t$, we get:
\begin{equation}
\sum_{t=1}^T \Big( f_t^c(x_t^0) - f_t^c(\hat{x}_t^\star) \Big) \le \sum_{t=1}^T \Big( \hat{f}_t^c(x_t^0) - \hat{f}_t^c(\hat{x}_t^\star) \Big) + T\delta^2\beta d. \label{eqn:app C7}
\end{equation}

{\it Step 4: Analyze the OCO Update on the Smoothed Cost.}
We now analyze the core OCO regret term $\sum_{t=1}^T \big( \hat{f}_t^c(x_t^0) - \hat{f}_t^c(\hat{x}_t^\star) \big)$ in the right-hand side of \eqref{eqn:app C7}.
We use the auxiliary function (which has been defined in {\it step 1.0}) $h_t(x)=\hat{f}_t^c(x)+(\tilde{g}_t-\nabla\hat{f}_t^c(x_t^0))x$, where $\tilde{g}_t=\frac{1}{2\delta}\big[f_t^c(\bar{x}_t^0)-f^c_t(\tilde{x}_t^0)\big]u_t^0$.
By construction, $\E_{u_t}[h_t(x)] = \hat{f}_t^c(x)$ for any $x$ independent of $u_t$. Thus, $\E[h_t(x_t^0)] = \E[\hat{f}_t^c(x_t^0)]$ and $\E[h_t(\hat{x}_t^\star)] = \E[\hat{f}_t^c(\hat{x}_t^\star)]$. We therefore need to bound $\E\Big[\sum_{t=1}^T \big( h_t(x_t^0) - h_t(\hat{x}_t^\star) \big)\Big]$.

The update rule is $x_{t+1}^0=\mathbf{\Pi}_{(1-\xi)\X}(x_t^0-\eta_t(\nabla h_t(x_t^0)+e_t))$, where $e_t = g_t^0 - \tilde{g}_t = e_t^m + e_t^c$.
Since $\hat{f}_t^c(\cdot)$ is $\mu$-strongly convex (Lemma~\ref{lemma:properties of hat f_t^c}), $h_t(\cdot)$ is also $\mu$-strongly convex. By definition of strong convexity:
\begin{equation}
h_t(x_t^0)-h_t(\hat{x}_t^\star)\le\nabla h_t(x_t^0)^{\top}(x_t^0-\hat{x}_t^\star)-\frac{\mu}{2}\norm{x_t^0-\hat{x}_t^\star}^2. \label{eqn:strong convex of h}
\end{equation}
From the update rule and Lemma~\ref{lemma:pythagorean} (since $\hat{x}_t^\star \in (1-\xi)\X$):
\begin{align*}
\norm{x_{t+1}^0 - \hat{x}_t^\star}^2 &= \norm{\mathbf{\Pi}_{(1-\xi)\X}(x_t^0-\eta_t(\nabla h_t(x_t^0)+e_t)) - \hat{x}_t^\star}^2 \\
&\le \norm{x_t^0-\eta_t(\nabla h_t(x_t^0)+e_t) - \hat{x}_t^\star}^2 \\
&= \norm{x_t^0 - \hat{x}_t^\star}^2 - 2\eta_t (x_t^0 - \hat{x}_t^\star)^{\top}(\nabla h_t(x_t^0)+e_t) + \eta_t^2 \norm{\nabla h_t(x_t^0)+e_t}^2.
\end{align*}
Rearranging this gives:
\begin{align}
(x_t^0-\hat{x}_t^\star)^{\top}\nabla h_t(x_t^0) \le \frac{\norm{x_t^0 - \hat{x}_t^\star}^2 - \norm{x_{t+1}^0 - \hat{x}_t^\star}^2}{2\eta_t} + (\hat{x}_t^\star-x_t^0)^{\top}e_t + \frac{\eta_t}{2}\norm{\nabla h_t(x_t^0)+e_t}^2. \label{eqn:upper bound on inner product}
\end{align}
Substitute \eqref{eqn:upper bound on inner product} into \eqref{eqn:strong convex of h}:
\begin{align}\nonumber
h_t(x_t^0)-h_t(\hat{x}_t^\star) &\le \frac{\norm{x_t^0 - \hat{x}_t^\star}^2 - \norm{x_{t+1}^0 - \hat{x}_t^\star}^2}{2\eta_t} + (\hat{x}_t^\star-x_t^0)^{\top}e_t \\\nonumber
&\quad + \frac{\eta_t}{2}\norm{\nabla h_t(x_t^0)+e_t}^2 - \frac{\mu}{2}\norm{x_t^0-\hat{x}_t^\star}^2 \\\nonumber
&\le \frac{\norm{x_t^0 - \hat{x}_t^\star}^2 - \norm{x_{t+1}^0 - \hat{x}_t^\star}^2}{2\eta_t} + (1-\xi)D\norm{e_t} \\
&\quad + \eta_t\norm{\nabla h_t(x_t^0)}^2 + \eta_t\norm{e_t}^2 - \frac{\mu}{2}\norm{x_t^0-\hat{x}_t^\star}^2, \label{eqn:h_t distance to optimal}
\end{align}
where we used $\norm{a+b}^2 \le 2\norm{a}^2 + 2\norm{b}^2$ and $\norm{\hat{x}_t^\star-x_t^0} \le (1-\xi)D$ from Remark~\ref{remark:feasible small set}.

Now we account for the dynamic comparator $\hat{x}_t^\star$. Use the fact that $\norm{a-b}^2 \ge \norm{a-c}^2 - 2\langle a-b, c-b \rangle \ge \norm{a-c}^2 - 2\norm{a-c}\norm{c-b}$:
\begin{align}\nonumber
\norm{x_{t+1}^0 - \hat{x}_t^\star}^2 &= \norm{x_{t+1}^0 - \hat{x}_{t+1}^\star + \hat{x}_{t+1}^\star - \hat{x}_t^\star}^2 \\
&\ge \norm{x_{t+1}^0 - \hat{x}_{t+1}^\star}^2 - 2\norm{x_{t+1}^0 - \hat{x}_{t+1}^\star}\norm{\hat{x}_{t+1}^\star - \hat{x}_t^\star} \nonumber \\
&\ge \norm{x_{t+1}^0 - \hat{x}_{t+1}^\star}^2 - 2(1-\xi)D\norm{\hat{x}_t^\star - \hat{x}_{t+1}^\star}. \label{eqn:lower bound on distance to optimal}
\end{align}
Substituting \eqref{eqn:lower bound on distance to optimal} into \eqref{eqn:h_t distance to optimal}: \\(specifically, $-\norm{x_{t+1}^0 - \hat{x}_t^\star}^2 \le -\norm{x_{t+1}^0 - \hat{x}_{t+1}^\star}^2 + 2(1-\xi)D\norm{\hat{x}_t^\star - \hat{x}_{t+1}^\star}$)
\begin{align*}
h_t(x_t^0)-h_t(\hat{x}_t^\star) &\le \frac{\norm{x_t^0 - \hat{x}_t^\star}^2 - \norm{x_{t+1}^0 - \hat{x}_{t+1}^\star}^2}{2\eta_t} + \frac{(1-\xi)D\norm{\hat{x}_t^\star-\hat{x}_{t+1}^\star}}{\eta_t} \\
&\quad + (1-\xi)D\norm{e_t} + \eta_t\norm{\nabla h_t(x_t^0)}^2 + \eta_t\norm{e_t}^2 - \frac{\mu}{2}\norm{x_t^0-\hat{x}_t^\star}^2.
\end{align*}
Now, we sum from $t=1$ to $T$. With the choice $\eta_t = 1/(t\mu)$, we have $\frac{1}{2\eta_t} - \frac{\mu}{2} = \frac{t\mu}{2} - \frac{\mu}{2} = \frac{(t-1)\mu}{2} = \frac{1}{2\eta_{t-1}}$.
The potential function terms form a telescoping sum that is non-positive:
\begin{align*}
\sum_{t=1}^T \Big( \frac{\norm{x_t^0 - \hat{x}_t^\star}^2}{2\eta_t} &- \frac{\norm{x_{t+1}^0 - \hat{x}_{t+1}^\star}^2}{2\eta_t} - \frac{\mu}{2}\norm{x_t^0-\hat{x}_t^\star}^2 \Big)\\ &= \sum_{t=1}^T \Big( (\frac{1}{2\eta_t} - \frac{\mu}{2})\norm{x_t^0-\hat{x}_t^\star}^2 - \frac{\norm{x_{t+1}^0 - \hat{x}_{t+1}^\star}^2}{2\eta_t} \Big) \\
&= \sum_{t=1}^T \Big( \frac{\norm{x_t^0-\hat{x}_t^\star}^2}{2\eta_{t-1}} - \frac{\norm{x_{t+1}^0 - \hat{x}_{t+1}^\star}^2}{2\eta_t} \Big) \\
&= \frac{\norm{x_1^0-x_{1,\xi}^{\star}}^2}{2\eta_0} - \frac{\norm{x_{T+1}^0 - \hat{x}_{t+1}^\star}^2}{2\eta_T} \le 0,
\end{align*}
where we set $1/\eta_0 = 0$ and $\hat{x}_{t+1}^\star=\hat{x}_t^\star$ since $t=T+1$ is vacuous in our analysis.

Thus, we arrive at the inequality:
\begin{align}
\sum_{t=1}^T\big(h_t(x_t^0)-h_t(\hat{x}_t^\star)\big)&\le\sum_{t=1}^{T-1}\frac{(1-\xi)D\norm{\hat{x}_t^\star-\hat{x}_{t+1}^\star}}{\eta_t}+\sum_{t=1}^T(1-\xi)D\norm{e_t} \nonumber \\
& \quad +\sum_{t=1}^T\eta_t\norm{\nabla h_t(x_t^0)}^2+\sum_{t=1}^T\eta_t\norm{e_t}^2 .\label{eqn:sum h_t distance to optimal solution}
\end{align}

{\it Step 5: Bound the Terms in the Inequality \eqref{eqn:sum h_t distance to optimal solution}.}
We now take the expectation of \eqref{eqn:sum h_t distance to optimal solution} and use the bounds derived in Appendix~\ref{app:omitted proofs in Thm 1} (which we cite here):
\begin{itemize}
    \item $\E\Big[\sum_{t=1}^{T-1}\frac{(1-\xi)D\norm{\hat{x}_t^\star-\hat{x}_{t+1}^\star}}{\eta_t}\Big] \le \frac{\sqrt{2}(1-\xi)D}{\mu}\E\Big[\sum_{t=1}^{T-1}\norm{\hat{x}_t^\star-\hat{x}_{t+1}^\star}\Big] = \frac{\sqrt{2}(1-\xi)^2D}{\mu}V_{T}$ (from \eqref{eqn:upper bound on y_t-y_t+1}).
    \item $\E\Big[\sum_{t=1}^T\eta_t\norm{\nabla h_t(x_t^0)}^2\Big] \le \frac{hG^2}{2(2h-1)\mu}(1+\ln T)$ (from \eqref{eqn:sum upper bound term 1}).
    \item $\E\Big[\sum_{t=1}^T(1-\xi)D\norm{e_t}\Big] \le (\frac{(1-\xi)D}{\delta(2(2h-1))^{1/4}}+\frac{\sqrt{2}h^2G(1-\xi)D}{2\delta(2(2h-1))^{1/2}})\sum_{t=1}^T\varphi_t+\frac{h^3G^2(1-\xi)D}{\delta(2(2h-1))^{1/2}}(1+\ln{T})$ (from \eqref{eqn:sum upper bound term 2}).
    \item $\E\Big[\sum_{t=1}^T\eta_t\norm{e_t}^2\Big] \le \frac{2}{\delta^2\mu\sqrt{2h-1}}\sum_{t=1}^T\varphi_t^2+\frac{4G^2h^2}{\mu(2h-1)}(1+\ln T)$ (from \eqref{eqn:sum upper bound term 3}).
\end{itemize}

Summing these bounds gives the bound for $\E\Big[\sum_{t=1}^T \big( \hat{f}_t^c(x_t^0) - \hat{f}_t^c(\hat{x}_t^\star) \big)\Big]$ (note that $\E[h_t(x_t^0)] = \E[\hat{f}_t^c(x_t^0)]$ and $\E[h_t(\hat{x}_t^\star)] = \E[\hat{f}_t^c(\hat{x}_t^\star)]$ which have been argued before):
\begin{align}\nonumber
\E\Big[\sum_{t=1}^T \big( \hat{f}_t^c(x_t^0) &- \hat{f}_t^c(\hat{x}_t^\star) \big)\Big]\\\nonumber
&\le \frac{\sqrt{2}(1-\xi)^2D}{\mu}V_{T} + \Big(\frac{hG^2}{2(2h-1)\mu} + \frac{h^3G^2(1-\xi)D}{\delta(2(2h-1))^{1/2}} + \frac{4G^2h^2}{\mu(2h-1)}\Big)(1+\ln T) \\\nonumber
& \quad + \frac{2}{\delta^2\mu\sqrt{2h-1}}\sum_{t=1}^T\varphi_t^2 + \Big(\frac{(1-\xi)D}{\delta(2(2h-1))^{1/4}}+\frac{\sqrt{2}h^2G(1-\xi)D}{2\delta(2(2h-1))^{1/2}}\Big)\sum_{t=1}^T\varphi_t \\\nonumber
&\le \frac{\sqrt{2}(1-\xi)^2D}{\mu}V_{T} + \Big(\frac{8G^2h^2+hG^2}{2\mu(2h-1)}+\frac{h^3G^2(1-\xi)D}{\delta(2(2h-1))^{1/2}}\Big)(1+\ln T) \\
& \quad + \frac{2}{\delta^2\mu\sqrt{2h-1}}\sum_{t=1}^T\varphi_t^2 + \Big(\frac{(1-\xi)D}{\delta(2(2h-1))^{1/4}}+\frac{\sqrt{2}h^2G(1-\xi)D}{2\delta(2(2h-1))^{1/2}}\Big)\sum_{t=1}^T\varphi_t.\label{eqn:sum f_t distance to optimal solution0}
\end{align}

{\it Combine All Pieces for Algorithmic Regret.}
We now combine the results from {\it Steps 2, 3, and 5}:
\begin{align*}
\E\Big[ C_T(\bm{x}^{0}) - C_T(\hat{\bm{x}}^\star) \Big] &\le \E\Big[\sum_{t=1}^T \Big( f_t^c(x_t^0) - f_t^c(\hat{x}_t^\star) \Big)\Big] + \E\Big[\Big|\sum_{t=1}^Tf_t(x_{t-h+1:t}^0)-f_t^c(x_t^0)\Big|\Big]\\&\qquad\qquad+ \E\Big[\Big|\sum_{t=1}^T f_t(\hat{x}_{t-h+1:t}^\star)-f_t^c(\hat{x}_t^\star)\Big|\Big] \\
&\overset{\eqref{eqn:app C7}}\le \E\Big[\sum_{t=1}^T \Big( \hat{f}_t^c(x_t^0) - \hat{f}_t^c(\hat{x}_t^\star) \Big)\Big] + T\delta^2\beta d + \E\Big[\Big|\sum_{t=1}^Tf_t(x_{t-h+1:t}^0)-f_t^c(x_t^0)\Big|\Big]\\&\qquad\qquad+ \E\Big[\Big|\sum_{t=1}^T f_t(\hat{x}_{t-h+1:t}^\star)-f_t^c(\hat{x}_t^\star)\Big|\Big] .
\end{align*}
Plugging everything in:
\begin{align*}
\E\Big[ C_T(\bm{x}^{0}) - C_T(\hat{\bm{x}}^\star) \Big] &\le \Big[\text{Bound from \eqref{eqn:sum f_t distance to optimal solution0}}\Big] + T\delta^2\beta d + \Big[\text{Bound from \eqref{eqn:app C8}}\Big]\\& + \Big[\text{Bound from \eqref{eqn:app C9}}\Big] \\
&\le \Big[ \frac{\sqrt{2}(1-\xi)^2D}{\mu}V_{T} + \Big(\frac{8G^2h^2+hG^2}{2\mu(2h-1)}+\frac{h^3G^2(1-\xi)D}{\delta(2(2h-1))^{1/2}}\Big)(1+\ln T) \\
& \quad + \frac{2}{\delta^2\mu\sqrt{2h-1}}\sum_{t=1}^T\varphi_t^2 + \Big(\frac{(1-\xi)D}{\delta(2(2h-1))^{1/4}}+\frac{\sqrt{2}h^2G(1-\xi)D}{2\delta(2(2h-1))^{1/2}}\Big)\sum_{t=1}^T\varphi_t \Big] \\
& \quad + T\delta^2\beta d \\
& \quad + \Big[ \frac{3G^2h^3}{\mu(2(2h-1))^{1/2}}(1+\ln T)+\frac{\sqrt{2}Gh^2}{\delta\mu(2(2h-1))^{1/4}}\sum_{t=1}^T\varphi_t \Big] \\
& \quad + (1-\xi)Gh^2 V_{T}.
\end{align*}
Rearranging terms gives the bound for the Algorithmic Regret:
\begin{align}\nonumber
\E\Big[ C_T(\bm{x}^{0}) &- C_T(\hat{\bm{x}}^\star) \Big]\le \Big(\frac{\sqrt{2}(1-\xi)^2D}{\mu}+(1-\xi)Gh^2\Big)V_{T} + T\delta^2\beta d \\\nonumber
& \quad + \Big(\frac{8G^2h^2+hG^2}{2\mu(2h-1)}+\frac{h^3G^2(1-\xi)D}{\delta(2(2h-1))^{1/2}}+\frac{3G^2h^3}{\mu(2(2h-1))^{1/2}}\Big)(1+\ln T) \\\nonumber
& \quad + \frac{2}{\delta^2\mu\sqrt{2h-1}}\sum_{t=1}^T\varphi_t^2 \\\nonumber
& \quad + \Big(\frac{(1-\xi)D}{\delta(2(2h-1))^{1/4}}+\frac{\sqrt{2}h^2G(1-\xi)D}{2\delta(2(2h-1))^{1/2}}+\frac{\sqrt{2}Gh^2}{\delta\mu(2(2h-1))^{1/4}}\Big)\sum_{t=1}^T\varphi_t\\\nonumber 
&\le \Big(\frac{\sqrt{2}(1-\xi)D}{\mu}+Gh^2\Big)V_{T,\xi} + T\delta^2\beta d  \\\nonumber
& \quad + \Big(\frac{8G^2h^2+hG^2}{2\mu(2h-1)}+\frac{h^3G^2(1-\xi)D}{\delta(2(2h-1))^{1/2}}+\frac{3G^2h^3}{\mu(2(2h-1))^{1/2}}\Big)(1+\ln T) \\ 
& \quad + \frac{2}{\delta^2\mu\sqrt{2h-1}}\sum_{t=1}^T\varphi_t^2 + \Big(\frac{\sqrt{2}h^2G(1-\xi)D}{2\delta(2(2h-1))^{1/2}}+\frac{\sqrt{2}Gh^2+(1-\xi)D\mu}{\delta\mu(2(2h-1))^{1/4}}\Big)\sum_{t=1}^T\varphi_t \label{eqn:final alg regret}.
\end{align}

{\bf Part 2: Bounding the Approximation Error.}
This part bounds $C_T(\hat{\bm{x}}^\star) - C_T(\bm{x}^{\star})$.
Since $C_T(\cdot)$ is $G\sqrt{Th}$-Lipschitz (from Lemma~\ref{lemma:properties of hat f_t^c}(b)):
\begin{align}\nonumber
|C_T(\hat{\bm{x}}^\star) - C_T(\bm{x}^{\star})| &\le G\sqrt{Th} \norm{\hat{\bm{x}}^\star - \bm{x}^{\star}} \\\nonumber
&= G\sqrt{Th} \sqrt{\sum_{t=1}^T \norm{(1-\xi)x_t^\star - x_t^\star}^2} \\\nonumber
&= G\sqrt{Th} \sqrt{\sum_{t=1}^T \xi^2 \norm{x_t^\star}^2} \\\nonumber
&\le G\sqrt{Th} \cdot \xi \sqrt{T D^2} \quad (\text{since } \norm{x_t^\star} \le D \text{ as } x_t^\star \in \mathcal{X}) \\
&= \xi G T \sqrt{h} D. \label{eqn:final approx error}
\end{align}

{\bf Final Bound.} Combining the bounds for Part 1 (Eq. \eqref{eqn:final alg regret}) and Part 2 (Eq. \eqref{eqn:final approx error}), we get the final bound for the total initialization regret $\E[R(\A_\text{init})]$:
\begin{align*}
\E[R(\A_\text{init})] &= \E\Big[ C_T(\bm{x}^{0}) - C_T(\hat{\bm{x}}^\star) \Big] + \Big( C_T(\hat{\bm{x}}^\star) - C_T(\bm{x}^{\star}) \Big) \\
&\le  \Big(\frac{\sqrt{2}(1-\xi)^2D}{\mu}+(1-\xi)Gh^2\Big)V_{T} + T\delta^2\beta d \\
& \quad + \Big(\frac{8G^2h^2+hG^2}{2\mu(2h-1)}+\frac{h^3G^2(1-\xi)D}{\delta(2(2h-1))^{1/2}}+\frac{3G^2h^3}{\mu(2(2h-1))^{1/2}}\Big)(1+\ln T) \\
& \quad + \frac{2}{\delta^2\mu\sqrt{2h-1}}\sum_{t=1}^T\varphi_t^2  + \Big(\frac{\sqrt{2}h^2G(1-\xi)D}{2\delta(2(2h-1))^{1/2}}+\frac{\sqrt{2}Gh^2+(1-\xi)D\mu}{\delta\mu(2(2h-1))^{1/4}}\Big)\sum_{t=1}^T\varphi_t  \\
& \quad + \xi G T \sqrt{h} D,
\end{align*}
which completes the proof of Theorem~\ref{thm:initial regret bound}.

\subsection{Proof of Theorem~\ref{thm:overall convergence}}

For our analysis in this proof, denote
\begin{equation}\label{eqn:def of tilde g}
\widetilde{\bm{g}}^j_s=\sum_{k=s}^{s+h-1}\frac{f_k(\bar{x}^j_{k-h+1:k})-f_k(\tilde{x}^j_{k-h+1:k})}{2\delta^{\prime}}u_s^j.
\end{equation}
We then have from Lemma~\ref{lemma:zeroth-order gradient} that $\E_{u_s^j}[\widetilde{\bm{g}}_s^j]=\partial\hat{C}_T(\bm{x}^j)/\partial x_s$ for all $s\in[T]$ and all $j\in\{0,\dots,K-1\}$ and thus $\nabla\hat{C}_T(\bm{x}^j)=\E_{\bm{u}^j}[\widetilde{\bm{g}}^j]$, where $\bm{u}^j=(u_1^j,\dots,u_T^j)\sim\widetilde{\N}(0,I_{dT})$ and $\widetilde{\bm{g}}^j=(\widetilde{\bm{g}}_1^j,\dots,\widetilde{\bm{g}}_T^j)$.

Let $\bm{x}_{\xi^\prime}^{\star}\in\argmin_{\bm{x}\in((1-\xi^\prime)\X)^T}C_T(\bm{x})$. To proceed, we first establish the convergence to the proxy optimal solution $\bm{x}_{\xi^\prime}^{\star}$. Define the potential function:
\begin{equation}\label{eqn:def of W_j}
W_j=(1+\gamma)^j(C_T(\bm{x}^j)-C_T(\bm{x}_{\xi^\prime}^{\star})).
\end{equation}
We then have
\begin{align}\nonumber
W_{j+1}-W_j&=(1+\gamma)^j\big((\gamma+1)(C_T(\bm{x}^{j+1})-C_T(\bm{x}^j))+\gamma(C_T(\bm{x}^j)-C_T(\bm{x}_{\xi^\prime}^{\star}))\big)\\
&=(1+\gamma)^j\big(C_T(\bm{x}^{j+1})-C_T(\bm{x}^j)+\gamma(C_T(\bm{x}^{j+1})-C_T(\bm{x}_{\xi^\prime}^{\star}))\big).\label{eqn:potential difference}
\end{align}
To upper bound $W_{j+1}-W_j$, we first prove in Appendix~\ref{app:omitted proofs in Thm2} that the following inequality holds for all $j\in\{0,\dots,K-1\}$ and all $\bm{y}\in((1-\xi^\prime)\X)^T$:
\begin{multline}
C_T(\bm{x}^{j+1})-C_T(\bm{y})\le\beta^{\prime}\langle\bm{x}^j-\bm{x}^{j+1}, \bm{x}^j-\bm{y}\rangle-\frac{\beta^{\prime}}{2}\norm{\bm{x}^{j+1}-\bm{x}^j}^2-\frac{\mu}{2}\norm{\bm{x}^j-\bm{y}}^2\\+\langle\nabla C_T(\bm{x}^j)-\bm{g}^j,\bm{x}^{j+1}-\bm{y}\rangle.\label{eqn:relation based on hat C_T}
\end{multline}
Letting $\bm{y}=\bm{x}^j$ (which holds for $j\ge 1$, and we analyze $j=0$ separately if needed, though the proof structure allows it) and $\bm{y}=\bm{x}_{\xi^\prime}^{\star}$ in \eqref{eqn:relation based on hat C_T} (note that $\bm{x}_{\xi^\prime}^{\star} \in ((1-\xi^\prime)\X)^T$), we respectively obtain the following:
\begin{align}
C_T(\bm{x}^{j+1})-C_T(\bm{x}^{j})&\le-\frac{\beta^{\prime}}{2}\norm{\bm{x}^{j+1}-\bm{x}^j}^2+\langle\nabla C_T(\bm{x}^j)-\bm{g}^j,\bm{x}^{j+1}-\bm{x}^j\rangle,\label{eqn:y=x_t}\\\nonumber
C_T(\bm{x}^{j+1})-C_T(\bm{x}_{\xi^\prime}^{\star})&\le\beta^{\prime}\langle\bm{x}^j-\bm{x}^{j+1},\bm{x}^j-\bm{x}_{\xi^\prime}^{\star}\rangle-\frac{\beta^{\prime}}{2}\norm{\bm{x}^{j+1}-\bm{x}^j}^2-\frac{\mu}{2}\norm{\bm{x}^{j}-\bm{x}_{\xi^\prime}^{\star}}^2\\
&\qquad\qquad\qquad\qquad\qquad\qquad+\langle\nabla C_T(\bm{x}^j)-\bm{g}^j,\bm{x}^{j+1}-\bm{x}_{\xi^\prime}^{\star}\rangle.\label{eqn:y=x^star}
\end{align}
Moreover, for any $\bm{y}=(y_1,\dots,y_T)$ with $y_t\in(1-\xi^\prime)\X$ for all $t\in[T]$, we have the following relation proved in Appendix~\ref{app:omitted proofs in Thm2}:
\begin{equation}\label{eqn:upper bound on epsilon}
\E\big[\langle\nabla C_T(\bm{x}^j)-\bm{g}^j,\bm{x}^{j+1}-\bm{y}\rangle\big]\le\varepsilon.
\end{equation}
Now, substituting \eqref{eqn:y=x_t}-\eqref{eqn:y=x^star} with the upper bound $\varepsilon$ into \eqref{eqn:potential difference}, we show in Appendix~\ref{app:omitted proofs in Thm2} that
\begin{equation}
\label{eqn:upper bound on potential difference}
\E_{\bm{u}^j}[W_{j+1}-W_j]\le (1+\gamma)^{j+1}\varepsilon.
\end{equation}
From \eqref{eqn:def of W_j} and \eqref{eqn:upper bound on potential difference}, we get
\begin{align*}
&\E_{\bm{u}^j}\big[C_T(\bm{x}^{j+1})-C_T(\bm{x}_{\xi^\prime}^{\star})\big]\le\frac{1}{1+\gamma}\big(C_T(\bm{x}^{j})-C_T(\bm{x}_{\xi^\prime}^{\star})\big)+\varepsilon\quad\forall j\in\{0,\dots,K-1\},\\
\Rightarrow & \E_{\CU_{K-1}}\big[C_T(\bm{x}^K)-C_T(\bm{x}_{\xi^\prime}^{\star})\big]\le\Big(\frac{1}{1+\gamma}\Big)^K\big(C_T(\bm{x}^0)-C_T(\bm{x}_{\xi^\prime}^{\star})\big)+\sum_{j=0}^{K-1}\Big(\frac{1}{1+\gamma}\Big)^{j}\varepsilon\\
\Rightarrow & \E_{\CU_{K-1}}\big[C_T(\bm{x}^K)-C_T(\bm{x}_{\xi^\prime}^{\star})\big]\le\Big(\frac{1}{1+\gamma}\Big)^K\E\big[C_T(\bm{x}^0)-C_T(\bm{x}_{\xi^\prime}^{\star})\big]+\frac{\varepsilon}{\gamma},
\end{align*}
where $\CU_{K-1}\triangleq(\bm{u}^0,\dots,\bm{u}^{K-1})$ encapsulates all the randomness of Algorithm~\ref{alg:two point feedback}.

\textbf{Overall Bound on the Regret $R(\A)$.} We now combine this result with the initialization regret to bound the total dynamic regret against the true optimal $\bm{x}^\star$. We decompose the regret as:
\begin{align*}
\E[R(\A)] = \E\big[C_T(\bm{x}^K) - C_T(\bm{x}^\star)\big] = \underbrace{\E\big[C_T(\bm{x}^K) - C_T(\bm{x}_{\xi^\prime}^{\star})\big]}_{\text{Term A}} + \underbrace{\Big(C_T(\bm{x}_{\xi^\prime}^{\star}) - C_T(\bm{x}^\star)\Big)}_{\text{Term B}}.
\end{align*}
\textit{Bounding Term B (Approximation Error):}
Since $\bm{x}^\star$ is in $\X^T$ and contains the origin, $(1-\xi^\prime)\bm{x}^\star \in ((1-\xi^\prime)\X)^T$. The approximation error is bounded by:
\begin{align*}
    \text{Term B} &= C_T(\bm{x}_{\xi^\prime}^{\star}) - C_T(\bm{x}^\star) \\&\overset{\text{(a)}}\le C_T((1-\xi^\prime)\bm{x}^\star) - C_T(\bm{x}^\star) \\&\overset{\text{(b)}}\le \xi^\prime G T \sqrt{h} D,
\end{align*}
where (a) is due to $\bm{x}_{\xi^\prime}^{\star}\in\argmin_{\bm{x}\in((1-\xi^\prime)\X)^T}C_T(\bm{x})$; (b) follows the same proof logic as \eqref{eqn:final approx error}.

\noindent\textit{Bounding Term A and Combining:}
We substitute the bound for Term A derived before:
$$
\text{Term A} \le \Big(\frac{1}{1+\gamma}\Big)^K \E\big[C_T(\bm{x}^0) - C_T(\bm{x}_{\xi^\prime}^{\star})\big] + \frac{\varepsilon}{\gamma}.
$$
We bound the initial error term $\E\big[C_T(\bm{x}^0) - C_T(\bm{x}_{\xi^\prime}^{\star})\big]$ by relating it to the total initialization regret. Since $C_T(\bm{x}_{\xi^\prime}^{\star}) \ge C_T(\bm{x}^\star)$, we have:
\begin{align*}
\E\big[C_T(\bm{x}^0) - C_T(\bm{x}_{\xi^\prime}^{\star})\big] &= \E\big[C_T(\bm{x}^0) - C_T(\bm{x}^\star)\big] - \big(C_T(\bm{x}_{\xi^\prime}^{\star}) - C_T(\bm{x}^\star)\big) \\
&\le \E[R(\A_\text{init})].
\end{align*}
Substituting these back into the regret decomposition, we obtain the final result:
\begin{align*}
\E[R(\A)] \le \Big(\frac{1}{1+\gamma}\Big)^K \E[R(\A_\text{init})] + \frac{\varepsilon}{\gamma} + \xi^\prime G T \sqrt{h} D.
\end{align*}
This completes the proof of Theorem~\ref{thm:overall convergence}.

\section{Omitted Details for the Proof of Theorem~\ref{thm:initial regret bound}}\label{app:omitted proofs in Thm 1}
\subsection{Proof of the bounds of the four terms in \eqref{eqn:sum h_t distance to optimal solution}} 
First, note that 
\begin{align}\nonumber
&\sum_{t=1}^{T-1}\frac{(1-\xi)D\norm{\hat{x}_t^\star-\hat{x}_{t+1}^\star}}{\eta_t}\overset{(a)}\le\frac{(1-\xi)D}{\mu}\sqrt{(\sum_{t=1}^{T-1}\norm{\hat{x}_t^\star-\hat{x}_{t+1}^\star}^2)(\sum_{t=1}^{T-1}\frac{1}{t^2})}\\\nonumber
\overset{(b)}\le&\frac{\sqrt{2}(1-\xi)D}{\mu}\sqrt{\sum_{t=1}^{T-1}\norm{\hat{x}_t^\star-\hat{x}_{t+1}^\star}^2}\\
\le&\frac{\sqrt{2}(1-\xi)D}{\mu}\sum_{t=1}^{T-1}\norm{\hat{x}_t^\star-\hat{x}_{t+1}^\star}=\frac{\sqrt{2}(1-\xi)^2D}{\mu}V_T,\label{eqn:upper bound on y_t-y_t+1}
\end{align}
where (a) follows from the CS-inequality and (b) follows from the fact
\begin{align}
\sum_{t=1}^{T-1}\frac{1}{t^2}\le 1+\sum_{t=2}^{T-1}\frac{1}{t(t-1)}=1+\sum_{t=2}^{T-1}(\frac{1}{t-1}-\frac{1}{t})\le2.\label{eqn:sum of 1/t square}
\end{align}

Next, we have from \eqref{eqn:upper bound on tilde g_t} that 
\begin{align}\nonumber
\sum_{t=1}^T\eta_t\norm{\nabla h_t(x_t^0)}^2&\le \frac{hG^2}{\mu}\sum_{t=1}^T\frac{1}{t}\norm{u_t^0}^4\\\nonumber&\overset{(a)}\le\frac{hG^2}{2(2h-1)\mu}\sum_{t=1}^T\frac{1}{t}\\
&\overset{(b)}
\le\frac{hG^2}{2(2h-1)\mu}(1+\ln T),\label{eqn:sum upper bound term 1}
\end{align}
where (a) follows from the fact that $\lVert u_t^0\rVert\le1/(2(2h-1))^{1/4}$, and (b) uses the fact $\sum_{t=1}^T\frac{1}{t}<1+\ln T$.

To proceed, we upper bound $\norm{e_t}$ for all $t\in[T]$, where $e_t=e_t^m+e_t^c$. We first have  from \eqref{eqn:def of e_t^m} that 
\begin{equation}\label{eqn:upper bound on e_t^m}
\norm{e_t^m}\le\frac{\norm{u^0}}{\delta}\varphi_t\le\frac{\varphi_t}{\delta(2(2h-1))^{1/4}}.
\end{equation}
To upper bound $\norm{e_t^c}$, we have from \eqref{eqn:def of e_t^c} that 
\begin{equation}\label{eqn:overall upper bound on e_t^c}
\norm{e_t^c}\le\frac{1}{2\delta}|f_t(\bar{x}_{t-h+1:t}^0)-f_t^c(\bar{x}_t^0)|\norm{u^0}+\frac{1}{2\delta}|f_t(\tilde{x}_{t-h+1:t}^0)-f_t^c(\tilde{x}_t^0)|\norm{u^0}.
\end{equation}
We may further upper bound bound $\norm{e_t^c}$ as
\begin{align}\nonumber
\norm{e_t^c}&\le\frac{1}{2\delta}|f_t(\bar{x}_{t-h+1:t}^0)-f_t(\tilde{x}_{t-h+1:t}^0)|\norm{u^0}+\frac{1}{2\delta}|f_t^c(\bar{x}_{t}^0)-f_t^c(\tilde{x}_{t}^0)|\norm{u^0}\\\nonumber
&\overset{(a)}\le\frac{G}{2\delta}\sum_{j=0}^{h-1}\norm{\bar{x}_{t-j}^0-\tilde{x}_{t-j}^0}\norm{u^0}+\frac{Gh}{2\delta}\norm{\bar{x}_t^0-\tilde{x}_t^0}\norm{u^0}\\\nonumber
&\overset{(b)}\le G\sum_{j=0}^{h-1} \norm{u^0}\norm{u^0}+Gh\norm{u^0}^2\\
&\le\frac{2Gh}{\sqrt{2(2h-1)}},\label{eqn:upper bound on e_t^c 2nd}
\end{align}
where (a) again follows from Assumption~\ref{ass:objective functions}(c), and (b) uses the relation between $\bar{x}_t^0$ and $\tilde{x}_t^0$. The above inequality further implies 
\begin{align}
\norm{e_t^c}^2\le\frac{2G^2h^2}{(2h-1)}.\label{eqn:upper bound on e_t^c square}
\end{align}

We then upper bound $\sum_{t=1}^T\eta_t\norm{e_t}^2$ in \eqref{eqn:sum h_t distance to optimal solution}. We have 
\begin{align}\nonumber
\sum_{t=1}^T\eta_t\norm{e_t}^2&\le\sum_{t=1}^T\eta_t(2\norm{e_t^m}^2+2\norm{e_t^c}^2)
\\\nonumber
&\overset{(a)}\le\sum_{t=1}^T\eta_t\frac{2\varphi_t^2}{\delta^2\sqrt{2(2h-1)}}+4\sum_{t=1}^T\eta_t\frac{G^2h^2}{(2h-1)}\\\nonumber
&\overset{(b)}\le\frac{2}{\delta^2\mu\sqrt{2(2h-1)}}\sqrt{(\sum_{t=1}^T\varphi_t^4)(\sum_{t=1}^T\frac{1}{t^2})}+\frac{4G^2h^2}{\mu(2h-1)}\sum_{t=1}^T\frac{1}{t}\\
&\overset{(c)}\le\frac{2}{\delta^2\mu\sqrt{2h-1}}\sum_{t=1}^T\varphi_t^2+\frac{4G^2h^2}{\mu(2h-1)}(1+\ln T),\label{eqn:sum upper bound term 3}
\end{align}
where (a) follows from \eqref{eqn:upper bound on e_t^m} and \eqref{eqn:upper bound on e_t^c square}, (b) follows from the choice $\eta=\frac{1}{t\mu}$ and the CS-inequality. To obtain (c), we used the fact $\sum_{i=1}^na_i^2\le(\sum_{i=1}^na_i)^2$ for $a_i\in\R^+$, $\sum_{t=1}^T1/t^2\le2$ from \eqref{eqn:sum of 1/t square} and $\sum_{t=1}^T\frac{1}{t}\le 1+\ln T$. Similarly, we can upper bound $\sum_{t=1}^T\eta_t\norm{e_t}$ as
\begin{align}\nonumber
\sum_{t=1}^T\eta_t\norm{e_t}&\le\sum_{t=1}^T\eta_t(\norm{e_t^m}+\norm{e_t^c})\\\nonumber
&\le \sum_{t=1}^{T} \frac{1}{t\mu } \frac{\varphi _t}{\delta (2(2h-1))^\frac{1}{4} } +\sum_{t=1}^{T} \frac{1}{t\mu } \frac{\sqrt{2}Gh}{\sqrt[]{2h-1} }\\
&\le\frac{\sqrt{2}}{\delta\mu(2(2h-1))^{1/4}}\sum_{t=1}^T\varphi_t+\frac{\sqrt{2}Gh}{\mu \sqrt{2h-1}}(1+\ln T).\label{eqn:upper bound on sum eta_t time e_t}
\end{align}

Finally, we upper bound $\sum_{t=1}^TD\norm{e_t}$ in \eqref{eqn:sum h_t distance to optimal solution}. We first note from Assumption~\ref{ass:objective functions}(c) that 
\begin{align*}
\quad|f_t(\bar{x}_{t-h+1:t}^0)-f_t^c(\bar{x}_t^0)|&\le G\sum_{j=1}^{h-1}\norm{\bar{x}_t^0-\bar{x}_{t-j}^0}\\
&\le G\sum_{j=1}^{h-1}\sum_{l=1}^j\norm{\bar{x}_{t-l+1}^0-\bar{x}_{t-l}^0}\\
&\le G\sum_{j=1}^{h-1}\sum_{l=1}^j\norm{x_{t-l+1}^0+\delta u^0-x_{t-l}^0-\delta u^0}\\
&\overset{(a)}\le G\sum_{j=1}^{h-1}\sum_{l=1}^j\big(\eta_{t-l}\norm{\nabla h_{t-l}(x_{t-l}^0)}+\eta_{t-l}\norm{e_{t-l}}\big)\\
&\overset{\eqref{eqn:upper bound on tilde g_t}}\le\sqrt{h}G^2\sum_{j=1}^{h-1}\sum_{l=1}^j\eta_{t-l}\norm{u^0}^2+G\sum_{j=1}^{h-1}\sum_{l=1}^j\eta_{t-l}\norm{e_{t-l}},
\end{align*}
where (a) uses the relation $x_{t-l+1}^0=\Pi_{(1-\xi)\X}(x_{t-l}^0-\eta_{t-l}
(\nabla h_{t-l}(x_{t-l}^0)+e_t))$ and Lemma~\ref{lemma:pythagorean} in Appendix~\ref{app:tech lemmas}. 

Similarly, one can show that 
\begin{align*}
|f_t(\tilde{x}_{t-h+1:t}^0)-f_t^c(\tilde{x}_t^0)|\le\sqrt{h}G^2\sum_{j=1}^{h-1}\sum_{l=1}^j\eta_{t-l}\norm{u^0}^2+G\sum_{j=1}^{h-1}\sum_{l=1}^j\eta_{t-l}\norm{e_{t-l}}.
\end{align*}
We then have from \eqref{eqn:overall upper bound on e_t^c} that 
\begin{align}\nonumber
&\norm{e^c_t}\le \frac{\sqrt{h}G^2}{\delta}\sum_{j=1}^{h-1}\sum_{l=1}^j\eta_{t-l}\norm{u^0}^3+\frac{G}{\delta}\sum_{j=1}^{h-1}\sum_{l=1}^j\norm{e_{t-l}}\norm{u^0}\\\nonumber
&\quad \quad \le \frac{ \sqrt{h}h^2G^2}{2\delta } \sum_{t=1}^{T} \eta _t \norm{u^0}^3 +\frac{h^2G^2}{2\delta } \sum_{t=1}^{T} \eta _t \norm{e_t}\norm{u^0},\\\nonumber
\overset{\eqref{eqn:upper bound on sum eta_t time e_t}}\Rightarrow&\E\big[\norm{e_t^c}\big]\le \frac{\sqrt{h}h^2G^2}{2\delta\mu(2(2h-1))^{3/4}}(1+\ln{T})+\frac{h^2G}{2\delta(2(2h-1))^{1/4}}\big[\frac{\sqrt{2}}{\delta\mu(2(2h-1))^{1/4}}\sum_{t=1}^T\varphi_t \\\nonumber
&\quad\quad\quad\quad\quad\quad\quad\quad\quad\quad+\frac{\sqrt{2}Gh}{\mu \sqrt{2h-1}}(1+\ln T)\big]\\
&\quad\quad\quad \le \frac{\sqrt{2}h^2G}{2\delta(2(2h-1))^{1/2}}\sum_{t=1}^T\varphi_t+\frac{h^3G^2}{\delta(2(2h-1))^{1/2}}(1+\ln{T}).\label{eqn:overall upper bound on e_t^c exp}
\end{align}
Combining \eqref{eqn:upper bound on e_t^m} and \eqref{eqn:overall upper bound on e_t^c exp} and summing from $t=1$ to $t=T$, we get 
\begin{align}\nonumber
&\quad\E\Big[\sum_{t=1}^TD\norm{e_t}\Big]\\
&\le (\frac{D}{\delta(2(2h-1))^{1/4}}+\frac{\sqrt{2}h^2GD}{2\delta(2(2h-1))^{1/2}})\sum_{t=1}^T\varphi_t+\frac{h^3G^2D}{\delta(2(2h-1))^{1/2}}(1+\ln{T}).\label{eqn:sum upper bound term 2} 
\end{align}
Combining the upper bounds in \eqref{eqn:upper bound on y_t-y_t+1}, \eqref{eqn:sum upper bound term 1}, \eqref{eqn:sum upper bound term 3} and \eqref{eqn:sum upper bound term 2}, noting that $\delta=1/\sqrt{T}\le1$ and recalling that $\E[h_t(x_t^0)]=\E[\hat{f}_t^c(x_t^0)]$ and $\E[h_t(\hat{x}_t^\star)]=\E[\hat{f}_t^c(\hat{x}_t^\star)]$, we obtain \eqref{eqn:sum f_t distance to optimal solution0} from \eqref{eqn:sum h_t distance to optimal solution}.$\hfill\blacksquare$

\section{Omitted Details in the Proof of Theorem~\ref{thm:overall convergence}}\label{app:omitted proofs in Thm2}
\subsection{Proof of \eqref{eqn:relation based on hat C_T}}
We first express $C_T(\bm{x}^{j+1})-C_T(\bm{y})=C_T(\bm{x}^{j+1})-C_T(\bm{x}^j)+C_T(\bm{x}^j)-C_T(\bm{y})$. From $\beta^{\prime}$-smoothness of $C_T(\cdot)$, we have
\begin{equation*}
C_T(\bm{x}^{j+1})-C_T(\bm{x}^j)\le\langle\nabla C_T(\bm{x}^j),\bm{x}^{j+1}-\bm{x}^j\rangle+\frac{\beta^{\prime}}{2}\norm{\bm{x}^{j+1}-\bm{x}^j}^2.
\end{equation*}
From the $\mu$-strongly convexity of $C_T(\cdot)$, we have
\begin{equation*}
C_T(\bm{x}^{j})-C_T(\bm{y})\le\langle \nabla C_T(\bm{x}^j), \bm{x}^{j}-\bm{y}\rangle-\frac{\mu}{2}\norm{\bm{x}^j-\bm{y}}^2.
\end{equation*}
Summing the above two inequalities yields
\begin{equation}
C_T(\bm{x}^{j+1})-C_T(\bm{y})\le\langle\nabla C_T(\bm{x}^j),\bm{x}^{j+1}-\bm{y}\rangle+\frac{\beta^{\prime}}{2}\norm{\bm{x}^{j+1}-\bm{x}^j}^2-\frac{\mu}{2}\norm{\bm{x}^j-\bm{y}}^2. \label{equ:3}
\end{equation}
As we argued in the main body, denoting $\bm{z}^{j+1}=\bm{x}^j-\alpha\bm{g}^j$ with $\alpha=1/\beta^{\prime}$, we get that
\begin{align}\nonumber
\langle \nabla C_T(\bm{x}^j),\bm{x}^{j+1}-\bm{y}\rangle&=\langle \bm{g}^j+\nabla C_T(\bm{x}^j)-\bm{g}^j,\bm{x}^{j+1}-\bm{y}\rangle\\\nonumber
&=\langle\bm{g}^j,\bm{x}^{j+1}-\bm{y}\rangle+\langle\nabla C_T(\bm{x}^j)-\bm{g}^j,\bm{x}^{j+1}-\bm{y}\rangle\\\nonumber
&=\beta^{\prime}\langle\bm{x}^j
-\bm{z}^{j+1},\bm{x}^{j+1}-\bm{y}\rangle+\langle\nabla C_T(\bm{x}^j)-\bm{g}^j,\bm{x}^{j+1}-\bm{y}\rangle\\\nonumber
&\le\beta^{\prime}\langle\bm{x}^{j}-\bm{x}^{j+1},\bm{x}^{j+1}-\bm{y}\rangle+\langle\nabla C_T(\bm{x}^j)-\bm{g}^j,\bm{x}^{j+1}-\bm{y}\rangle\\\nonumber
&=\beta^{\prime}\langle\bm{x}^{j}-\bm{x}^{j+1},\bm{x}^{j+1}-\bm{x}^{j}+\bm{x}^{j}-\bm{y}\rangle+\langle\nabla C_T(\bm{x}^j)-\bm{g}^j,\bm{x}^{j+1}-\bm{y}\rangle\\\nonumber
&=-\beta^{\prime}\norm{\bm{x}^{j+1}-\bm{x}^j}^2+\beta^{\prime}\langle\bm{x}^{j}-\bm{x}^{j+1},\bm{x}^{j}-\bm{y}\rangle+\langle\nabla C_T(\bm{x}^j)-\bm{g}^j,\bm{x}^{j+1}-\bm{y}\rangle,
\end{align}
where the inequality follows from Lemma~\ref{lemma:pythagorean}, noting that $\bm{x}^{j+1}=\mathbf{\Pi}_{(1-\xi^\prime)\X}(\bm{z}^{j+1})$ and assuming $\bm{y} \in ((1-\xi^\prime)\X)^T$. This assumption holds for $\bm{y}=\bm{x}_{\xi^\prime}^\star$ and $\bm{y}=\bm{x}^j$ for $j \ge 1$ (and for $j=0$, it requires $(1-\xi)\X \subseteq (1-\xi^\prime)\X$, which can be easily achieved). Plugging the above into~\eqref{equ:3}, we get~\eqref{eqn:relation based on hat C_T}.$\hfill\blacksquare$

\subsection{Proof of \eqref{eqn:upper bound on epsilon}}
For any $\bm{y}=(y_1,\dots,y_T)$ with $y_t\in(1-\xi^\prime)\X$ for all $t\in[T]$, we have
\begin{align*}
\E_{\bm{u}^j}\big[\langle\nabla C_T(\bm{x}^j)-\bm{g}^j,\bm{x}^{j+1}-\bm{y}\rangle\big]&=\E_{\bm{u}^j}\big[\langle\nabla C_T(\bm{x}^j)-\widetilde{\bm{g}}^j+\widetilde{\bm{g}}^j-\bm{g}^j,\bm{x}^{j+1}-\bm{y} \rangle\big],\\
&=\underbrace{\E_{\bm{u}^j}\big[\langle\nabla C_T(\bm{x}^j)-\widetilde{\bm{g}}^j,\bm{x}^{j+1}-\bm{y}\rangle\big]}_{\text{(i)}}+\underbrace{\E_{\bm{u}^j}\big[\langle\widetilde{\bm{g}}^j-\bm{g}^j,\bm{x}^{j+1}-\bm{y}\rangle\big]}_{\text{(ii)}}.
\end{align*}
\normalsize
We first upper bound (ii). Recalling the definition of $\bm{g}^j$ (resp., $\widetilde{\bm{g}}^j$) given in Algorithm~\ref{alg:two point feedback} (resp., \eqref{eqn:def of tilde g}), we get that
\begin{align*}
\E_{{u}^j_s}[\bm{g}^j_s-\widetilde{\bm{g}}^j_s]=\E_{u^j_s}\Big[\sum_{k=s}^{s+h-1}\frac{l_k(\bar{x}_{k-h+1:k})-f_k(\bar{x}_{k-h+1:k})+f_k(\tilde{x}_{k-h+1:k})-l_k(\tilde{x}_{k-h+1:k})}{2\delta^{\prime}}u^j_s\Big],
\end{align*}
for all $s\in[T]$. By Assumption~\ref{ass:model error} and the fact that $\Vert u_s^j\Vert\le1/2(2h-1))^{1/4}$, we obtain
\begin{align*}
&\E_{\bm{u}^j}\big[\lVert\bm{g}_s^j-\widetilde{\bm{g}}_s^j\rVert\big]\le\frac{\sum_{k=s}^{s+h-1}\varphi_k}{\delta^{\prime}(2(2h-1))^{1/4}}\ \forall s\in[T],\\
\Rightarrow&\E_{\bm{u}^j}\big[\lVert \bm{g}^j-\widetilde{\bm{g}}^j\rVert\big]\le\frac{\sum_{s=1}^T\sum_{k=s}^{s+h-1}\varphi_k}{\delta^{\prime}(2(2h-1))^{1/4}}\le\frac{h\sum_{t=1}^T\varphi_t}{\delta^{\prime}(2(2h-1))^{1/4}}.
\end{align*}
Thus, we can upper bound (ii) as
\begin{align*}
\text{(ii)}\le\E_{\bm{u}^j}\big[\Vert\bm{g}^j-\widetilde{\bm{g}}^j\Vert\Vert\bm{x}^{j+1}-\bm{y}\Vert\big]\le \frac{(1-\xi^\prime)Dh\sum_{t=1}^T\varphi_t}{\delta^{\prime}(2(2h-1))^{1/4}},
\end{align*}
where the second inequality follows from Remark~\ref{remark:feasible small set} (noting $\bm{x}^{j+1}, \bm{y} \in ((1-\xi^\prime)\X)^T$).

We then upper bound (i). Note that
\begin{align*}
&\E_{\bm{u}^j}\big[\langle\nabla C_T(\bm{x}^j)-\widetilde{\bm{g}}^j,\bm{x}^{j+1}-\bm{y}\rangle\big]\\
&\le\E_{\bm{u}^j}\big[\Vert\nabla C_T(\bm{x}^j)-\widetilde{\bm{g}}^j\Vert\Vert\bm{x}^{j+1}-\bm{y}\Vert\big]\\
&\overset{(a)}\le (1-\xi^\prime)D\E_{\bm{u}^j}\big[\Vert\nabla C_T(\bm{x}^j)-\widetilde{\bm{g}}^j\Vert\big]\\
&\le (1-\xi^\prime)D\E_{\bm{u}^j}\Big[\sqrt{\Vert\nabla C_T(\bm{x}^j)-\widetilde{\bm{g}}^j\Vert^2}\Big]\\
&\overset{(b)}\le (1-\xi^\prime)D\sqrt{\E_{\bm{u}^j}\big[\Vert\nabla C_T(\bm{x}^j)-\widetilde{\bm{g}}^j\Vert^2\big]}\\
\Rightarrow&\frac{1}{(1-\xi^\prime)^2 D^2}\big(\E_{\bm{u}^j}\big[\langle\nabla C_T(\bm{x}^j)-\widetilde{\bm{g}}^j,\bm{x}^{j+1}-\bm{y}\rangle\big]\big)^2\\
&\le \E_{\bm{u}^j}\Big[\norm{\nabla C_T(\bm{x}^j)}^2+\norm{\widetilde{\bm{g}}^j}^2-2\nabla C_T(\bm{x}^j)^{\top}\widetilde{\bm{g}}^j\Big]\\
&=\norm{\nabla C_T(\bm{x}^j)}^2-2\nabla C_T(\bm{x}^j)^{\top}\nabla\hat{C}_T(\bm{x}^j)+\E_{\bm{u}^j}\big[\norm{\widetilde{\bm{g}}^j}^2\big]\\
&=\nabla C_T(\bm{x}^j)^{\top}\big(-\nabla C_T(\bm{x}^j)+2\nabla C_T(\bm{x}^j)-2\nabla\hat{C}_T(\bm{x}^j)\big)+\E_{\bm{u}^j}\big[\norm{\widetilde{\bm{g}}^j}^2\big]\\
&=-\norm{\nabla C_T(\bm{x}^j)}^2+2\norm{\nabla C_T(\bm{x}^j)}\norm{\nabla C_T(\bm{x}^j)-\nabla\hat{C}_T(\bm{x}^j)}+\E_{\bm{u}^j}\big[\norm{\widetilde{\bm{g}}^j}^2\big]\\
&\overset{(c)}\le-\norm{\nabla C_T(\bm{x}^j)}^2+\E_{\bm{u}^j}\big[\norm{\widetilde{\bm{g}}^j}^2\big]+G\beta hT(Td+3)^{3/2}\delta^{\prime},
\end{align*}
where (a) follows from Remark~\ref{remark:feasible small set}, (b) follows from Jensen's inequality and the concavity of $\sqrt{\cdot}$ and (c) follows from Lemma~\ref{lemma:properties of hat f_t^c}(b) and Lemma~\ref{lemma:zeroth-order gradient}(b).\footnote{Note the well-known fact that if a differentiable function $f:\R^d\to\R$ is convex and $G$-Lipschitz, then $\Vert\nabla f(x)\Vert\le G$ for all $x\in\R^d$. Thus, we know from Lemma~\ref{lemma:properties of hat f_t^c}(b) that $\Vert \nabla C_T(\bm{x})\Vert\le G\sqrt{Th}$ for all $\bm{x}\in\R^{dT}$.} It remains to upper bound the term $-\Vert C_T(\bm{x}^j)\Vert^2+\E_{\bm{u}^j}\big[\Vert\widetilde{\bm{g}}^j\Vert^2\big]$. Noting that $\sum_{k=s}^{s+h-1}f_k(x_{k-h+1:k})$ may be viewed as a function of $x_{s-h+1},\dots,x_{s+h-1}$, denoted as $F_s(x_{s-h+1:s+h-1})$, we see from \eqref{eqn:def of tilde g} that $\widetilde{\bm{g}}^j_s$ can be equivalently written as
\begin{equation*}
\widetilde{\bm{g}}^j_s=\frac{1}{2\delta^{\prime}}\big(F_s(\bar{x}^j_{s-h+1:s+h-1})-F_s(\tilde{x}^j_{s-h+1:s+h-1})\big)u_s^j.
\end{equation*}
Moreover, using similar arguments to those in the proof of Lemma~\ref{lemma:properties of hat f_t^c}, one can show that $F_s(\cdot)$ is $\beta h$-smooth. Now, following similar arguments to those in the proof of \cite[Theorem~4]{nesterov2017random}, one can show that
\begin{align*}
\norm{\widetilde{\bm{g}}^j_s}^2&\le\frac{\delta^{\prime 2}}{8}\beta^2h^2\norm{u_s^j}^6+2\norm{u_s^j}^4\norm{\nabla F_s(x^j_{s-h+1:s+h-1})}^2\\
&\le\frac{\delta^{\prime 2}\beta^2h^2}{8(2(2h-1))^{3/2}}+\frac{1}{(2h-1)}\norm{\nabla F_s(x^j_{s-h+1:s+h-1})}^2,
\end{align*}
where the second inequality follows from the fact that $\Vert u_s^j\Vert\le1/(2(2h-1))^{1/4}$. One can also show that $\sum_{s=1}^T\Vert \nabla F_s(x^j_{s-h+1:s+h-1})\Vert^2\le(2h-1)\Vert\nabla C_T(\bm{x}^j)\Vert^2$. It follows that 
\begin{align*}
\norm{\widetilde{\bm{g}}^j}^2=\sum_{s=1}^T\norm{\widetilde{\bm{g}}^j_s}^2\le\frac{\delta^{\prime 2}\beta^2h^2T}{8(2(2h-1))^{3/2}}+\norm{\nabla C_T(\bm{x}^j)}^2.
\end{align*}
Thus, we can upper bound $\text{(i)}+\text{(ii)}$ as
\begin{align*}
\text{(i)}+&\text{(ii)}\le\sqrt{\frac{(1-\xi^\prime)^2 D^2\delta^{\prime 2}\beta^2h^2T}{8(2(2h-1))^{3/2}}+(1-\xi^\prime)^2 D^2Gh\beta T(Td+3)^{3/2}\delta^{\prime}}+\frac{(1-\xi^\prime)Dh\sum_{t=1}^T\varphi_t}{\delta^{\prime}(2(2h-1))^{1/4}}\\
&\le\frac{(1-\xi^\prime)D\delta^{\prime}\beta h\sqrt{T}}{2\sqrt{2}(2(2h-1))^{3/4}}+\sqrt{h\beta GT\delta^{\prime}}(1-\xi^\prime)D(Td+3)^{3/4}+\frac{(1-\xi^\prime)Dh\sum_{t=1}^T\varphi_t}{\delta^{\prime}(2(2h-1))^{1/4}}\triangleq\varepsilon.
\end{align*}
Therefore, we complete the proof of \eqref{eqn:upper bound on epsilon}.$\hfill\blacksquare$

\subsection{Proof of \eqref{eqn:upper bound on potential difference}}
We get from \eqref{eqn:y=x_t}-\eqref{eqn:y=x^star} and the upper bound $\varepsilon$ in \eqref{eqn:upper bound on epsilon} that
\begin{align*}
&\quad C_T(\bm{x}^{j+1})-C_T(\bm{x}^j)+\gamma(C_T(\bm{x}^{j+1})-C_T(\bm{x}_{\xi^\prime}^{\star}))\\
&\le-\frac{\beta^{\prime}}{2}(1+\gamma)\norm{\bm{x}^{j+1}-\bm{x}^j}+\gamma\beta^{\prime}\langle\bm{x}^{j}-\bm{x}^{j+1},\bm{x}^j-\bm{x}_{\xi^\prime}^{\star}\rangle-\frac{\mu\gamma}{2}\norm{\bm{x}^{j}-\bm{x}_{\xi^\prime}^{\star}}^2+(1+\gamma)\varepsilon\\
&\le-\frac{1}{\beta^{\prime}-\mu}\norm{\beta^{\prime}(\bm{x}^{j+1}-\bm{x}^j)-\mu(\bm{x}^j-\bm{x}_{\xi^\prime}^{\star})}^2+(1+\gamma)\varepsilon \le (1+\gamma)\varepsilon,
\end{align*}
where the last inequality follows from $\gamma\triangleq\mu/(\beta^{\prime}-\mu)$ and some simple algebra, completing the proof of \eqref{eqn:upper bound on potential difference}.$\hfill\blacksquare$

\section{Technical Lemma}\label{app:tech lemmas}
\begin{lemma}\citep[Proposition~2.2]{bansal2019potential}
\label{lemma:pythagorean}
Given a convex set $\X\subseteq\R^d$, let $a\in\X$, $b^{\prime}\in\R^d$ and $b=\mathbf{\Pi}_{\X}(b^{\prime})$. Then $\langle a-b,b^{\prime}-b\rangle\le0$ and $\norm{a-b}^2\le\norm{a-b^{\prime}}^2$.
\end{lemma}

\end{document}